\documentclass[11pt]{amsart}
\usepackage{amsmath, amsthm, amssymb}
\usepackage{graphicx}
\usepackage{hyperref}
\usepackage[usenames]{color}
\usepackage{srcltx} 
\usepackage{verbatim}
\usepackage{bbm}

\definecolor{dgreen}{cmyk}{1,0,1,0.2}

\newtheorem{thm}{Theorem}[section]
\newcommand{\bt}{\begin{thm}}
\newcommand{\et}{\end{thm}}

\newtheorem{conj}[thm]{Conjecture}

\newtheorem{ex}[thm]{Example}

\newtheorem{cor}[thm]{Corollary}   
\newcommand{\bc}{\begin{cor}}
\newcommand{\ec}{\end{cor}}

\newtheorem{lem}[thm]{Lemma}   
\newcommand{\bl}{\begin{lem}}
\newcommand{\el}{\end{lem}}

\newtheorem{prop}[thm]{Proposition}
\newcommand{\bp}{\begin{prop}}
\newcommand{\ep}{\end{prop}}

\newtheorem{defn}[thm]{Definition}
\newcommand{\bd}{\begin{defn}}    
\newcommand{\ed}{\end{defn}}

\newtheorem{rmrk}[thm]{Remark}   

\newcommand{\br}{\begin{rmrk}}
\newcommand{\er}{\end{rmrk}}

\newcommand{\distto}{\stackrel {\textrm{dist}}{\longrightarrow} }

\newcommand{\Fto}{\stackrel {\mathcal{F}}{\longrightarrow} }
\newcommand{\VFto}{\stackrel {\mathcal{VF}}{\longrightarrow} }
\newcommand{\VADBto}{\stackrel {\mathcal{VADB}}{\longrightarrow} }

\newcommand{\be}{\begin{equation}}

\newcommand{\ee}{\end{equation}}

\newcommand{\diam}{\operatorname{Diam}}

\newcommand{\Ricci}{\rm{Ricci}}
\newcommand{\Scal}{{\rm Scal}} 
\newcommand{\disjointunion}{\sqcup}

\newcommand{\Lip}{\operatorname{Lip}}

\newcommand{\mass}{{\mathbf M}}

\newcommand{\area}{\operatorname{Area}}
\newcommand{\vol}{\operatorname{Vol}}

\usepackage{color, comment}

\title[Geometric Stability of the Zero Mass Rigidity Theorem]{Geometric Stability of the \\ Schoen-Yau Zero Mass Rigidity Theorem}

\author{Christina Sormani}
\address{CUNYGC and Lehman College, New York}
\email{\href{mailto:sormanic@gmail.com}{sormanic@gmail.com}}
\thanks{Sormani's work was supported by a Scholar Incentive Leave from Lehman College of CUNY, and the Simons Center for Geometry and Physics at SUNY Stony Brook.}

\begin{document}


\begin{abstract}
 In 1979, Schoen and Yau proved their famous Positive Mass Theorem which is a combination of a comparison theorem: {\em a three dimensional asymptotically flat Riemannian manifold with nonnegative scalar curvature has nonnegative ADM mass}, and a rigidity theorem: {\em if such a manifold has zero ADM mass then it is isometric to Euclidean space}.   Here we review results and open questions on the geometric stability  of their zero mass rigidity theorem: {\em if such a manifold has almost zero mass, how close is its geometry to that of Euclidean space}?  We review the geometry of these spaces, examples of sequences of such spaces with mass approaching zero, and a variety of 
 geometric notions of convergence. 
 Although there has been much progress, it is still an open question (even in dimension three): exactly which geometric notion of convergence works best to capture the geometric stability of this famous rigidity theorem.

{\bf \textcolor{blue}{Dedicated to Professor Shing-Tung Yau whose brilliant work and incredible mentorship has deeply influenced the lives and careers of so many mathematicians.
Without his support, I would have taken a job in industry and never even published my doctoral dissertation. He was an amazing postdoctoral mentor, dedicating so much time and energy to his research team when we worked with him and then providing advice and guidance throughout our careers.  I am deeply grateful to him.}}

\end{abstract}

\maketitle

 \newpage

\tableofcontents

\newpage

\section{Introduction}\label{sect:intro}

In this paper, we view the Schoen-Yau Positive Mass Theorem of \cite{Schoen-Yau-positive-mass} as a comparison theorem about the geometry of three dimensional asymptotically flat manifolds with nonnegative scalar curvature.  We consider the rigidity result within their famous theorem, that we call the Schoen-Yau Zero Mass Rigidity Theorem.  We survey open questions, examples, and theorems regarding the conjectured geometric stability of this powerful rigidity theorem. 

This is not an exhaustive survey of the many beautiful applications and results related to the Positive Mass Theorem.   Instead, we focus on the geometry of sequences of manifolds whose masses approach zero and explore whether such sequences converge in a way which controls their geometry. We include a review of all results regarding the Gromov-Hausdorff (GH), metric measure (mm), and volume preserving intrinsic flat  ($\mathcal VF$) limits of regions in these asymptotically flat manifolds with nonnegative scalar curvature and discuss conjectures related to these notions of convergence.

Some of the most elegant theorems in Riemannian Geometry are comparison theorems which compare the geometry of a complete manifold, $M^3$, with nonnegative curvature to the geometry of Euclidean space, ${\mathbb E}^3$.  For nonnegative sectional curvature, there is the Toponogov Triangle Comparison Theorem which provides an inequality involving distances in triangles.  
For nonnegative Ricci curvature there is the  Bishop Volume Comparison Theorem comparing the volumes of balls, $B_p(r)\subset M^3$, with
balls in Euclidean space, $B_0(r)\subset {\mathbb E}^3$: $\vol(B_p(r))\le \omega_nr^n$.
See the texts of Bishop-Crittenden and Cheeger-Ebin for detailed proofs of these and other such comparison theorems \cite{Bishop-Crittenden}\cite{Cheeger-Ebin}.

For nonnegative scalar curvature, we can view the
Schoen-Yau Positive Mass Theorem  as a geometric comparison theorem about the ADM mass of an asymptotically flat manifold.  We state this comparison theorem here in a very geometric way using the Huisken-Yau foliation  by constant mean curvature (CMC) surfaces, $\Sigma(R)$ of area, $\area(\Sigma(R))=4\pi R^2$ \cite{Huisken-Yau-center}.  These are now known to be outer boundaries of unique isoperimetric regions, $\Omega(R)$ by the work of Chodosh-Eichmair-Shi-Yu \cite{CESY-isoperimetry} and others (cited within).

\begin{thm}[Schoen-Yau Mass Comparison Theorem]
\label{thm:comp}
If $M^3$ is asymptotically flat,
in the sense that there are isoperimetric
regions, $\Omega(R)\subset M$, for $R>0$ sufficiently large, and diffeomorphisms,
\be\label{eq:end}
\exists F: M\setminus \Omega(R) \subset M \to  {\mathbb E}^3
\setminus B_0(R)
\ee
with sufficient decay of partial derivatives
in each order:
\be\label{eq:vague}
F_{*}g \to g_{{\mathbb E}^3} \textrm{ as } R\to \infty,
\ee
then
\be
\Scal \ge 0 \implies m_{ADM}(M^3) \ge 0.
\ee
\end{thm}

We are deliberately vague about asymptotic flatness here in (\ref{eq:vague})
as the strength of asymptotic flatness 
varies in different papers and the reader should consult original source material for details.
Recall that the ADM mass is
the limit of the Hawking masses
\be \label{eq:ADM-mass}
m_{ADM}(M^3)= \lim_{R\to\infty} m_H(\Sigma_R),
\ee 
of the CMC surfaces, $\Sigma_R$, because they are round spheres at infinity. The Hawking mass depends on mean curvature,
$H_\Sigma$, and area,
\be\label{eq:Hawking-mass}
m_H(\Sigma)=\sqrt{\tfrac{\area(\Sigma)}{16\pi}}
\left(1 - \tfrac{1}{16\pi}\int_\Sigma H_{\Sigma}^2 \, d\sigma\right).
\ee
See Nerz's work characterizing asymptotic flatness using the CMC foliation \cite{Nerz-Geometric}.

We review the geometry of manifolds with nonnegative scalar curvature in Section~\ref{sect:geometry} with subsections on distances, areas and minimal surfaces, isoperimetric regions, and notions of mass.    Note it is important to keep in mind there are many notions of quasi-local mass and it is useful to check each satisfies the properties described by Liu-Yau in \cite{Liu-Yau-2003,Liu-Yau-quasi}.  In this paper we focus on the most geometric notion of mass requiring the least differentiability.  We discuss Huisken's isoperimetric mass and related work of Benatti-Fogagnolo, Cederbaum-Nerz, Corvino-Wu, Huang, Jauregui-Lee, Miao, Shi, and others. We briefly touch on Gromov's $\mu$-bubbles and refer to Lesourd-Unger-Yau's work for more details.

The mass comparison theorem stated above holds for a larger class of manifolds via the 
Penrose Inequality which states that,
 \be \label{eq:Penrose}
\Scal \ge 0
\implies 
m_{ADM}(M) \ge 
\sqrt{\tfrac{\area(\partial M)}{16\pi}}\ge 0.
\ee
This was proven first by Huisken-Ilmanen \cite{Huisken-Ilmanen} assuming $\partial M$ is a connected outward minimizing surface, using inverse mean curvature flow and Geroch monotonicity of Hawking mass.  By work of Bray in \cite{Bray-Penrose}, the Penrose Inequality holds when $M^3$ is in the following class of manifolds:

\begin{defn} \label{defn:class}
Let ${\mathcal M}$ be the class of asymptotically flat three dimensional Riemannian manifolds with nonnegative scalar curvature such that $M$ has no closed interior minimal surfaces and the boundary, $\partial M$, may be empty or have multiple connected components (called horizons) each of which is a closed minimal surface. 
\end{defn}

When $M$ satisfying the hypothesis of Theorem~\ref{thm:comp} is not in the class $\mathcal{M}$ of Definition~\ref{defn:class}, one can always cut out all closed minimal surfaces in $M$ and then the connected component, $M'$, containing the end
is in this class: $M'\in \mathcal{M}$.   In Section~\ref{sect:examples}, we will review many examples of such manifolds, discuss why the Penrose Inequality only holds for manifolds in this class $\mathcal{M}$, and
describe the cutting process in more detail.

Ideally all comparison theorems have corresponding rigidity theorems which state that, when the equality is obtained, then the Riemannian manifold is isometric to a specific comparison space.  For example, the Bishop Volume Comparison Rigidity Theorem states that if $M^m$ is a complete Riemannian manifold 
with nonnegative Ricci curvature that has 
Euclidean volume growth, $\lim_{r\to \infty}\vol(B_p(r))/r^m=\omega_m$, then $M^m$ is isometric to ${\mathbb E}^m$.

For nonnegative scalar curvature and ADM mass, the Schoen-Yau Positive Mass Comparison Theorem of \cite{Schoen-Yau-positive-mass}
includes the following rigidity statement
which holds for the full class of 
$M\in \mathcal{M}$ by the work of Bray \cite{Bray-Penrose}:

\begin{thm}[Schoen-Yau Zero Mass Rigidity Theorem]
\label{thm:rigid}
If $M^3\in \mathcal{M}$ of Definition~\ref{defn:class},
then
\be\label{eq:scalar-zero-mass}
\Scal \ge 0 \textrm{ and } m_{ADM}(M^3) = 0
\ee
implies that $M^3$ is isometric to ${\mathbb E}^3$.
\end{thm}

An isometry, $\Psi: X \to Y$, completely identifies the geometry of the two spaces.  By definition, an isometry $\Psi$ is a bijection 
that is {\em distance preserving}:
\be \label{eq:dist-pres}
d_{X}(p,q)=d_{Y}(\Psi(p),\Psi(q))
\qquad \forall p,q \in X.
\ee
As a consequence, an isometry also preserves volumes of regions,
\be \label{eq:dist-vol}
\vol_X(\Omega)=\vol_Y(\Psi(\Omega))\qquad
\forall \Omega\subset X,
\ee
preserves boundaries and areas of boundaries, 
\be \label{eq:bndry-pres}
\Psi(\partial \Omega)=\partial \Psi(\Omega)
\textrm{ and }
\area_X(\partial \Omega)=\area_Y(\Psi(\partial \Omega)),
\ee
where volume and area are defined using Hausdorff measures when $X$ and $Y$ are only metric spaces (see within for details).
An isometry thus also preserves isoperimetric regions:
\be \label{eq:isoper-pres}
\textrm{
{\em $\Omega$ is isoperimetric in $X$
$\iff$ $\Psi(\Omega)$ is isoperimetric in $Y$}.} 
\ee
With sufficient regularity on $\Sigma$ and $X$, $\Psi$ preserves mean curvature,
\be \label{eq:mean-pres}
H_{\Sigma\subset X}(p)=H_{\Psi(\Sigma)\subset Y}(\Psi(p)),
\ee
and thus preserves minimal surfaces and constant mean curvature surfaces,
\be \label{eq:min-pres}
\textrm{
{\em $\Sigma$ is a min surf in $X$
$\iff$ $\Psi(\Sigma)$ is a min surf in $Y$},}
\ee
\be \label{eq:CMC-pres}
\textrm{
{\em $\Sigma$ is a CMC surface in $X$
$ \iff $ $\Psi(\Sigma)$ is a CMC surface in $Y$},}
\ee
as well as preserving Hawking masses and ADM mass, 
\be \label{eq:Hawking-pres}
m_H(\Sigma\subset X)=m_H(\Psi(\Sigma)\subset Y)
\textrm{ and }
m_{ADM}(X)=m_{ADM}(Y).
\ee
There is much research ongoing as to how much regularity is required to control mean curvature, minimal surfaces,
masses, and scalar curvature.
For settings with less regularity one can use Huisken's isoperimetric mass.

Ideally a geometric rigidity theorem has a corresponding geometric stability theorem which states that when the equality is almost achieved then the geometry of the space is almost the same as the geometry of the comparison space in the rigidity theorem.  Equivalently, a geometric stability theorem might start with a sequence of manifolds satisfying the hypotheses of the original comparison theorem which approach the equality in the rigidity theorem, and conclude that this sequence converges in some geometric way to the geometry of the comparison space in the rigidity theorem.  The goal is to use a geometric notion of convergence which allows one to also obtain approximations of all the geometric equalities in (\ref{eq:dist-pres})-(\ref{eq:Hawking-pres}) even in settings where one cannot hope to find a bijection.

The first geometric stability theorems are the almost rigidity theorems involving volume and Ricci curvature appearing in papers of Colding \cite{Colding-volume} and Cheeger-Colding \cite{ChCo-almost-rigidity} using Gromov-Hausdorff Convergence \cite{Gromov-metric}.   Their simplest geometric stability theorem states that if a sequence of complete Riemannian manifolds, $M_j^m$, have
$\Ricci\ge 0$ and the 
$
\vol(B(p_j,r)) \to \omega_m r^m
$
then 
the balls converge in the Gromov-Hausdorff (GH) sense to Euclidean balls.
By an example of Perelman in
\cite{Perelman-max-vol}, it is known these balls need not even have bijections to the Euclidean ball.   Colding also proved metric measure convergence in the sense of Fukaya \cite{Colding-volume}\cite{Fukaya-collapsing}.   Matveev-Portegies showed volume preserving intrinsic flat $(\mathcal{VF})$ convergence in the sense of Sormani-Wenger \cite{Matveev-Portegies}\cite{SW-JDG}.   We will review these geometric notions of convergence and their properties in Section~\ref{sect:distances}.  

There has been significant progress towards proving a Geometric Stability Theorem for the Zero Mass Rigidity Theorem.  We first state this conjecture vaguely here using the class of manifolds, $\mathcal{M}$ of Definition~\ref{defn:class}, from the Penrose Inequality.  For other Geometric Stability Conjectures involving scalar curvature see papers by Huisken-Ilmanen, Gromov, and the IAS survey on Scalar curvature \cite{Huisken-Ilmanen}\cite{Gromov-Plateau}\cite{Gromov-four}\cite{Sormani-IAS-survey}.

\begin{conj}[Geometric Stability of the Zero Mass Rigidity Theorem] \label{conj:geom}
If $M_j^3\in \mathcal{M}$
then
\be\label{eq:scalar-mass-to-zer}
\Scal \ge 0 \textrm{ and } m_{ADM}(M_j^3) \to 0
\ee
implies 
$
M^3_j \to  {\mathbb E}^3
$
with respect to some geometric notion of convergence.
\end{conj}

\begin{figure}[h] 
   \centering
\includegraphics[width=5in]{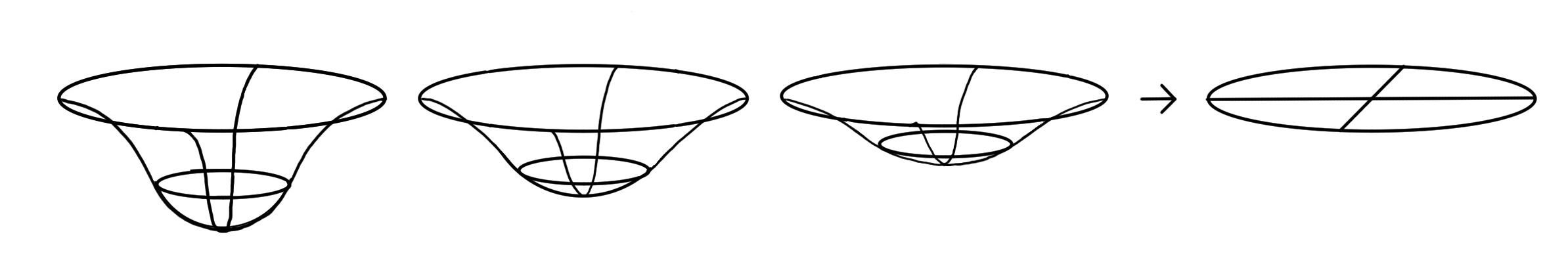} 
   \caption{A simple example of a sequence, $M_j^3$, of spherically symmetric asymptotically flat Riemannian manifolds with nonnegative scalar curvature such that $m_{ADM}(M_j^3) \to 0$ which smoothly converge to flat Euclidean space.
   }
\label{fig:smooth-limit}
\end{figure}

In Section~\ref{sect:examples} we present many examples of sequences of Riemannian manifolds which satisfy the hypotheses of this conjecture.  These include examples which converge smoothly like the one depicted in Figure~\ref{fig:smooth-limit} and far more complicated examples.  We begin with scalar flat examples including Riemannian Schwarzschild space
and spaces with many black holes first studied by Misner, Brill-Lindquist, and Lichnerowicz,
with more recent work by Stavrov et al (cited within).   We include well known examples with bubbling behind minimal surfaces that motivated the assumption that the manifolds lie in the class $\mathcal{M}$ of Definition~\ref{defn:class} that satisfies the Penrose Inequality.   There are also examples from Lee-Sormani's study of this conjecture in the spherically symmetric setting
demonstrating that this conjecture is false for Gromov-Hausdorff convergence.   We also present examples constructed using the Schoen-Yau tunnels of \cite{Schoen-Yau-tunnels} and sewing techniques by Basilio, Dodziuk, Sormani (cited within).   We ask a question about the scrunching of regions without tunnels.
We present the examples before defining notions of convergence so that readers may keep their minds open to the idea of developing a new notion of geometric convergence themselves which fits these examples.

In Section~\ref{sect:distances} we review various geometric notions of convergence for sequences of Riemannian manifolds including Gromov-Hausdorff (GH), metric measure (mm), intrinsic flat ($\mathcal{F}$), volume preserving intrinsic flat ($\mathcal{VF}$).  We define each notion and apply it to each sequence from Section~\ref{sect:examples} and then discuss the properties of each notion of convergence.  All of the notions above are geometric, so we contrast them with smoother notions of convergence where the entire sequence is diffeomorphic
starting with $C^k$, $C^0$, and VADB convergence.  
We point to literature contrasting Lebesgue and Sobolev convergence of metric tensors with the geometric notions of convergence.     We discuss Lee-Naber-Neumayer's suggested approach to abandon geometry and use $d_p$ distances rather than the geometric distances defined by lengths of shortest curves \cite{LNN-dp}.  We close with a discussion of convergence away from removed sets, including work of Lakzian-Sormani \cite{Lakzian-Sormani}, Dong \cite{Dong-some-stability} and  Dong-Song \cite{Dong-Song-Stability}.

Before proceeding to the paper, let us close the introduction with a precise restatement of the geometric stability conjecture that should be true but remains an open question:   

\begin{conj}[$\mathcal VF$ Stability of the Schoen-Yau Zero Mass Rigidity Theorem] \label{conj:LS}
If $M_j^3\in \mathcal{M}$ of Definition~\ref{defn:class} satisfy the hypotheses of
Conjecture~\ref{conj:geom}
with a uniform radius $R>0$ such that the isoperimetric regions
$\Omega_j(R)\subset M_j$ have uniformly bounded diameter,
\be
\diam_{g_j}(\Omega_j(R))\le D_R, 
\ee
and $M_j$ have ends that
are uniformly converging to
the Euclidean end,
\be\label{eq:endj}
\exists F_j: M_j\setminus \Omega_j(R) \subset M_j \to  r^{-1}[R,\infty)\subset {\mathbb E}^3
\ee
such that
\be
|F_{j*}g_j-g_{{\mathbb E}^3}|_{C^2(r^{-1}[R,\infty))}\to 0,
\ee
then the isoperimetric regions converge
in the volume preserving intrinsic flat sense to Euclidean balls,
\be
\Omega_j(R) \VFto B(0,R).
\ee
\end{conj}

This conjecture was first stated in this way in the final  section of Lee-Sormani \cite{LeeSormani1} along with a number of other related conjectures.  It was proven there in the spherically symmetric setting.  
As there have been so many advances in the understanding of isoperimetric sets and their relationship with mass and with intrinsic flat convergence, it is time for us to reconsider this version of the geometric stability of the zero mass conjecture. 

Throughout this paper we discuss why we have stated the conjecture in this way and why we believe it is true as stated.   We discuss partial results towards this conjecture and counter examples that support the choice of the hypotheses.  We also discuss how one might find a counter example and how one might design a different notion of geometric convergence.  We try to review all work towards this particular statement of the geometric stability conjecture.  However, we do not discuss all the many exciting approaches being taken towards proving Conjecture~\ref{conj:geom}.

We are quite deliberately avoiding stating precise definitions and theorems and precise conditions on examples in this survey as we are providing an intuitive discussion.  All readers should go to the original source material for precise statements of theorems and refer to definitions within the source material which may vary significantly from paper to paper.   This is particularly important when referring to asymptotic flatness, hypotheses on regularity, and notions of mass and of convergence.

\vspace{.3cm}
\noindent
{\bf Acknowledgements:}  I gratefully acknowledge support from the Simons Center for Geometry and Physics, Stony Brook University, at which some of the research for this paper was performed when I was on a Scholar Incentive Leave from CUNY.  This survey arises out of many discussions there in the the Fall of 2025.    I would especially like to thank Jeff Jauregui and Raquel Perales for their detailed written description of the literature that they shared with me after the SCGP program.  I would also like to thank Allen, Jauregui, Lee and Perales for feedback on an earlier version of this survey.  

\section{The Geometry of Manifolds with Nonnegative Scalar Curvature}\label{sect:geometry}

In this paper we consider complete asymptotically flat three-dimensional manifolds with nonnegative scalar curvature.
The precise definition of asymptotic flatness varies from paper to paper, so we suggest the reader refer to the original sources to verify the necessary asymptotic flatness conditions required for each result.  As we are primarily focused upon the Positive Mass Theorem as a geometric comparison theorem, we will focus on the geometric quantities that are involved in various statements of this theorem. 
Additional geometric quantities controlled by nonnegative scalar curvature, even in a compact setting, are reviewed in the IAS Scalar Curvature Survey \cite{Sormani-IAS-survey} and other sources mentioned there.

In Section~\ref{sect:distance} we quickly review the properties of distances on Riemannian manifolds.
In Section~\ref{sect:minimal} we first review the Schoen-Yau Theorem on stable closed minimal surfaces and the concept of outer minimizing surfaces.  In Section~\ref{sect:isoper} we review the Huisken-Yau CMC foliation of the asymptotically flat end and newer results on isoperimetric regions.  In Section~\ref{sect:mass} we briefly review various notions of mass and quasi-local mass.
We recommend reading Dan Lee's textbook \cite{Lee-text} for a more thorough review of these notions.  

\subsection{Distances}\label{sect:distance} 

Lets provide a brief guide to metric geometry aimed at the trained Riemannian geometer, as we will need metric geometry later to define geometric notions of convergence.
Given a Riemannian manifold, $M$, with a Riemannian metric tensor, $g$, we define the Riemannian distance between points, $p,q\in M$ by taking the infimum of lengths of pieceweise $C^1$ curves between the points:
\be \label{eq:Riem-dist-1}
d_{g}(p,q)=\inf \{L_{g}(C)\,|\, C:[a,b] \to M, \,\,C(a)=p,\, C(b)=q\}
\ee
where
\be \label{eq:Riem-dist-2}
L_{g}(C)=\int_a^b |C'(s)|_{g} \, ds.
\ee
In this way we convert the Riemannian manifold, $(M,g)$ into a metric space, $(M,d_g)$.   If $M$ is compact, then so is this metric space.

Given a metric space, $(M,d)$, we can define the ``rectifiable length'' of a continuous curve,
$C:[a,b]\to M$ as a supremum,
\be
L_d(C) =\sup\sum_{i=1}^{N}
d(C(t_i),C(t_{i-1}))
\ee
where the supremum is taken over all partitions
\be
a=t_0<t_1<\cdots <t_N=b.
\ee
We say a curve is rectifiable if it has finite length.
Taking the infimum over all such rectifiable lengths of rectifiable curves defines the ``intrinsic length distance'' on the metric space,
\be\label{eq:length-dist}
d_{L_d}(p,q)=
\inf \{L_{d}(C)\,|\, C:[a,b] \to M, \,\,C(a)=p,\, C(b)=q\}
\ee
which is not necessarily the original distance, $d$.  See Example~\ref{ex:sphere}. 

\begin{ex}\label{ex:sphere}
Consider for example: the metric space defined by taking the standard sphere, ${\mathbb S}^2$ with the Euclidean distance, $d=d_{{\mathbb E}^3}$, between the points. The rectifiable lengths of smooth curves in the sphere agree with the usual lengths
of the curves defined using the standard 
metric tensor, $g_{{\mathbb S}^2}$.  The 
intrinsic length distance between points is
achieved by geodesic arcs so that it agrees with the Riemannian distance on the sphere,
\be
d_{L_{d_{{\mathbb E}^3}}}=d_{g_{{\mathbb S}^2}}.
\ee
\end{ex}

On a path connected compact metric space, the infimum in (\ref{eq:length-dist}) is achieved by a ``length minimizing curve'' if 
\be
L_d(C)=d_{L_d}(C(a),C(b))
\ee
which is parametrized by arclength
\be
L_d(C([0,s])=d(C(0),C(s)) \quad \forall s\in [a,b].
\ee
Its important to note that metric geometers will sometimes call such curves ``geodesics" but that is problematic for Riemannian geometers, because in Riemannian geometry a geodesic is only locally length minimizing.  Recall that a geodesic in Riemannian geometry is a smooth curve with zero acceleration.  

If we consider
$\Omega\subset M$,  we should keep in mind that $\Omega$ has an intrinsic distance
\be\label{eq:int-dist}
d_{g,\Omega}(p,q)=
\inf \{L_{g}(C)\,|\, C:[a,b] \to \Omega, \,\,C(a)=p,\, C(b)=q\}
\ee
and an extrinsic distance
\be\label{eq:ext-dist}
d_{g,M}(p,q)=
\inf \{L_{g}(C)\,|\, C:[a,b] \to M, \,\,C(a)=p,\, C(b)=q\}.
\ee
See Example~\ref{ex:sphere} where
we contrast these intrinsic and extrinsic distances for ${\mathbb S}^2 \subset {\mathbb E}^3$.   

In this paper we are concerned with a
compact Riemannian manifold with boundary,
$(\Omega,g)$, which we will
convert into a compact metric space, $(\Omega, d_{g,\Omega})$. If the manifold has a convex boundary, then the length minimizing geodesics between points in the manifold are indeed geodesics of zero acceleration of the Riemannian manifold.   However if the boundary is not convex, part of the length minimizing curve can crawl along this boundary.   The curve can be broken into components, some of which are geodesics in the interior of $\Omega$ 
with zero acceleration and some of which are geodesics in the boundary,
$\partial \Omega$, with acceleration perpendicular to the boundary.

Whenever one is studying geometric notions one must be careful with the word isometric embedding.   In Riemannian geometry, an isometric embedding, $F:M\to N$,
is a smooth map such that the metric tensor pushes forward and so lengths of curves are preserved.  However the distances need not be preserved.  In Example~\ref{ex:sphere}, we saw that the Riemannian isometric embedding, $\iota:{\mathbb S}^2\to {\mathbb E}^3$, does not preserve distances.  In fact,
\be
d_{{\mathbb E}^3}(\iota(p),\iota(q))
< d_{{\mathbb S}^2}(p,q).
\ee

It is a beautiful theorem that a bijection, $F:M\to N$, between Riemannian manifolds is a Riemannian isometry, $F_*g_M=g_N$, if and only if it is distance preserving
\be
d_{g_N}(F(p),F(q))=d_{g_M}(p,q) \qquad \forall p,q \in M.
\ee
For more about Riemannian manifolds and metric spaces we refer the reader to the metric geometry text by Burago-Burago-Ivanov \cite{BBI}.

Also of importance is the notion of a Lipschitz map between two metric spaces,
$F:(M,d_M) \to (N,d_N)$ which is a map with
finite Lipschitz constant,
\be\label{eq:Lip}
\Lip(F)=\sup_{p\neq q \in M} \frac{d_N(F(p),F(q))}{d_M(p,q)}<\infty.
\ee
A bi-Lipschitz map is a bijection whose inverse is also Lipschitz.  

Once one has defined distances on a metric space, one can define the Hausdorff measures of a set, $U$, inside the metric space, $(M,d)$ as follows:
\be \label{eq:Hausdorff-measure}
{\mathcal{H}}^n_d(U) 
= \lim_{\delta\to 0}  \inf\left\{c_n \sum_{j=1}^\infty (\diam(U_j))^n\, : \, U \subset \bigcup_{j=1}^\infty U_j, \,\, \diam(U_j) \le \delta\right\}
\ee
which can be infinite or zero.  So a smooth surface, $\Sigma^2\subset M^3$, has 
Hausdorff measure,
$
{\mathcal{H}}_{d_g}^2(\Sigma)=\area_g(\Sigma)
$
where the area is defined as usual using integration of the metric tensor. 

In a metric space, a set $U$ has Hausdorff dimension,
\be \label{eq:Hausdorff-dim}
dim_{\mathcal H}(U)=\inf\{n\ge 0\,:\,
{\mathcal{H}}_d^n(U)=0\}.
\ee
So a smooth surface, $\Sigma^2\subset M^3$, has Hausdorff dimension $dim_{\mathcal H}(\Sigma^2)=2$.

 Distances are strongly controlled under nonnegative sectional curvature bounds as seen in the Toponogov Triangle Comparison Theorem.  This leads to the notion of an Alexandrov Space in metric geometry.   For nonnegative Ricci curvature, distances are weakly controlled because volumes of balls and annuli are controlled in the Bishop and Bishop-Gromov Volume Comparison Theorems and these balls and annuli involve distances.    

\begin{rmrk}\label{rmrk:surfaces}
The only theorems about distances in manifolds with nonnegative scalar curvature involve surfaces. Schoen-Yau have theorems about the distances and minimal surfaces in regions with $\Scal\ge \Lambda>0$ in \cite{Schoen-Yau-black}.   Gromov-Lawson \cite{Gromov-Lawson-positive}, Marques-Neves
\cite{Marques-Neves-rigidity}, and Liokumovich-Maximo \cite{Liokumovich-Maximo-waist} bound the diameters of surfaces in fibrations of such a manifold.  See also Remark~\ref{rmrk-mu-bubbles}.
\end{rmrk}

\subsection{Areas and Minimal Surfaces}
\label{sect:minimal}

The first geometric theorem about scalar curvature is the Schoen-Yau Stable Minimal Surface Theorem of
\cite{Schoen-Yau-minimal}.  Recall that a closed minimal surface, $\Sigma^2\subset M^3$,
is a compact surface without boundary that
is locally area-minimizing. That is,
around any point on $\Sigma^2$, there is a sufficiently small region $U\subset M$, 
such that
\be
\area(\Sigma^2\cap U) \le \area(\Sigma') \qquad \forall \Sigma'\subset U \,\,s.t.\,\, \partial \Sigma'=\partial(\Sigma^2\cap U).
\ee
This is well defined for 2 dimensional surfaces, $\Sigma$, on metric spaces using the second Hausdorff measure and Hausdorff dimension as in (\ref{eq:Hausdorff-measure})-(\ref{eq:Hausdorff-dim}).   On a smooth Riemannian manifold, minimal surfaces have zero mean curvature, 
\be
H_{\Sigma\subset M}(p)=0 \qquad \forall p \in \Sigma.
\ee

A closed minimal surface $\Sigma\subset M$ is stable if for any continuous deformation, $\Sigma_t$, 
through
$\Sigma_0=\Sigma^2$,
 there is an $\epsilon>0$,
such that 
\be
\area(\Sigma^2) \le \area(\Sigma^2_t) \qquad \forall t\in (-\epsilon, \epsilon).
\ee
Again this is well defined on metric spaces
 using the second Hausdorff measure and Hausdorff dimension as in (\ref{eq:Hausdorff-measure})-(\ref{eq:Hausdorff-dim}).

In \cite{Schoen-Yau-minimal}, Richard Schoen and Shing-Tung Yau proved that if $g$ is $C^2$ then by taking the second variation of the area and applying the Gauss-Bonnet Theorem, the only stable closed minimal surfaces in a manifold with nonnegative scalar curvature are tori and spheres.  
The Schoen-Yau proof of their positive mass theorem  involves the construction of a sequence of stable minimal surfaces with boundary
and a limiting complete minimal surface \cite{Schoen-Yau-positive-mass}. 
 See also work of Fischer-Colbrie and Schoen on complete minimal surfaces \cite{Fischer-Colbrie-Schoen}.   
 These original articles are beautiful and worth reading.

The Penrose Inequality bounds the mass of an asymptotically flat manifold with nonnegative scalar curvature
in terms of the area of an outermost minimal surface called an apparent horizon
as in (\ref{eq:Penrose}).   An outermost minimal surface, $\Sigma_{end}$, is a closed minimal surface that depends on the choice of a complete asymptotically flat end in $M$. It is defined by the property that the connected component $M'_{end}$ of $M\setminus \Sigma_{end}$
that contains the given end, has no closed minimal surfaces.  Generally the end is assumed to be given.  The Penrose Inequality is false if one does not assume the surface is outward minimizing as can be seen in examples with bubbles as depicted in Figure~\ref{fig:one-bubble}
which has a larger minimal surface of large area hidden behind the outermost minimal surface. 
See Remark~\ref{rmrk:Penrose-bubble}.  These outermost minimal surfaces are viewed as the apparent horizons of black holes in Penrose's intuitive derivation of the inequality.  Huisken-Ilmanen \cite{Huisken-Ilmanen} and Bray \cite{Bray-Penrose} prove the Penrose Inequality rigorously with geometric arguments.

\subsection{Isoperimetric regions}
\label{sect:isoper}

There are a number of beautiful theorems about isoperimetric regions and constant mean curvature for asymptotically flat manifolds with nonnegative scalar curvature.
Recall that an isoperimetric region, $\Omega^3(R)\subset M^3$, is a region whose outer boundary area is the area of a sphere of radius $R$ in Euclidean space, ${\mathbb E}^3$:
\be \label{eq:OmegaR}
\area (\partial \Omega^3(R)\setminus \partial M^3)=4\pi R^2
\ee
and whose volume is the maximum volume over all $\Omega$ satisfying this area constraint.
This is well defined on metric spaces
 using Hausdorff measures  as in (\ref{eq:Hausdorff-measure})-(\ref{eq:Hausdorff-dim}) in the place of volumes and areas.  
 
 The outer boundary, $\Sigma(R)=\partial \Omega(R)\setminus \partial M$, is called an isoperimetric surface.  On a smooth Riemannian manifold, a smooth isoperimetric surface has constant mean curvature (denoted CMC).     Huisken-Yau proved that on an asymptotically flat manifold with positive scalar curvature, for large enough $R$, there is a unique foliation by CMC surfaces that foliate the end \cite{Huisken-Yau-center}.  These CMC surfaces can then be used to show the ADM mass is defined canonically as in the introduction (\ref{eq:ADM-mass})-(\ref{eq:Hawking-mass}) and also to define a center of mass.  
 See also the work of Huang \cite{Huang-Foliations},
 Corvino-Wu \cite{Corvino-Wu}, 
 and of Cederbaum-Nerz \cite{Cederbaum-Nerz}  which involve various rates of asymptotic flatness.   Nerz characterizes asymptotic flatness using the CMC foliation in \cite{Nerz-Geometric}.
 
Bray proved  that the Hawking mass is
monotone increasing along the isoperimetric constant mean curvature
surfaces and converges to the ADM mass at infinity \cite{Bray-Penrose}.  This can be contrasted with the Geroch monotonicity of Hawking masses of the level sets of Huisken-Ilmanen's inverse mean curvature flow that we discussed in the introduction \cite{Geroch-monotonicity} \cite{Huisken-Ilmanen}.

In \cite{Shi-isoper}, Shi proved if $R$ is sufficiently large so that the outer boundary is in the asymptotically flat region, then 
\be
\vol(\Omega(R))\ge 4\pi R^3/3
\ee
with equality iff $\Omega(R)$ is isometric to the Euclidean ball, $B_0(R)$.  See also a new proof in a paper by Cheng \cite{Cheng-power}.  

Applying this result and the Huisken-Ilmanen inverse mean curvature flow,
Chodosh, Eichmair, Shi, and Yu prove in 
\cite{CESY-isoperimetry} that the canonical CMC foliation of the end by Huisken-Yau in \cite{Huisken-Yau-center} are the boundaries of isoperimetric regions whenever the manifold has nonnegative scalar curvature and positive mass. In fact they are uniquely isoperimetric for their given areas.   Antonelli-Fogagnolo-Nardulli-Pozzeta lower the regularity required to achieve these results in 
\cite{AFNP-positive}. See also the survey of Benatti-Fogagnolo  \cite{BF-isoper-survey}.

Note that this result and that of Shi in \cite{Shi-isoper} are results for large isoperimetric regions with outer boundaries in the asymptotically flat end.   It would be interesting to explore the properties of smaller isoperimetric regions and the isoperimetric profile.   See work of Nardulli and Osario Acevedo for small areas \cite{Nardulli-OsarioAcevedo}.  See work of
Antonelli-Fogagnolo-Nardulli-Pozzetta in \cite{AFNP-positive} and of Xu in \cite{Xu-isoperimetry} on the isoperimetric profiles of asymptotically flat manifolds with nonnegative scalar curvature.   

While all the results stated above are published using isoperimetric regions defined for fixed volume and minimizing area, we fix areas and maximize the volume  because of examples constructed in Lee-Sormani \cite{LeeSormani1} and included in Section~\ref{sect:examples} below where the boundaries of fixed areas are converging nicely rather than the interiors.

\begin{rmrk}\label{rmrk-mu-bubbles}
Gromov introduced the notion of a $\mu$ bubble (which is a prescribed mean curvature surface) in \cite{Gromov-positive}.  Distances have been estimated using these surfaces in work of Gromov \cite{Gromov-metric-ineq}\cite{Gromov-four}, Zhu \cite{Zhu-rigidity}\cite{Zhu-width},
and Lesourd-Unger-Yau
\cite{Lesourd-Unger-Yau}.   Distances involving capillary $\mu$ bubbles
appear in work of Wu \cite{Wu-capillary} building on work of Li \cite{Li-poly}.
\end{rmrk}

\begin{rmrk}\label{rmrk-FF}
It is worth recalling that the field of Geometric Measure Theory is being applied when one studies isoperimetric domains.
Federer and Fleming first applied their notion of integral currents and flat convergence to prove existence of isoperimetric sets in \cite{FF}. 
 The notion of integral currents was adapted to metric spaces by Ambrosio-Kirchheim in \cite{AK}.
They proved the same compactness theorem in their setting and also studied isoperimetric sets in metric spaces.   Sormani-Wenger intrinsic flat convergence is defined building upon these great works in \cite{SW-JDG}.
\end{rmrk}

\subsection{Notions of Mass}
\label{sect:mass}

The original definition of the ADM mass of an asymptotically flat manifold with nonnegative scalar curvature by Arnowitt-Desner-Misnor in \cite{ADM-mass} is defined using a coordinate chart and then proven to be independent of the choice of chart.  This is the definition used by Schoen-Yau and Witten in their original proofs of the Positive Mass Theorem \cite{Schoen-Yau-positive-mass}\cite{Witten-positive-mass}.

Once Huisken and Yau introduced the unique foliation
of an asymptotically flat end by stable CMC surfaces in \cite{Huisken-Yau}, one could see the more geometric
formula for ADM mass using Hawking masses
that we reviewed in
(\ref{eq:ADM-mass})-(\ref{eq:Hawking-mass}). 
 This coordinate-free approach depends only on the areas and mean curvatures of these surfaces.   The only concern is that 
 computing the mean curvature requires more regularity than we might have available to us when studying limits of sequences of manifolds. 

Huisken's isoperimetric mass is defined using only the isoperimetric deficits of regions:
\be
m_{ISO}(\Omega)=\frac{2}{\area(\partial \Omega)} \left( \vol(\Omega)-\frac{\area(\partial \Omega)^{3/2}}{6\sqrt{\pi}} \right).
\ee
Following the work of
Chodosh, Eichmair, Shi, and Yu in \cite{CESY-isoperimetry}, we can write the definition of this mass using their isoperimetric regions as 
in (\ref{eq:OmegaR}) as follows: 
\be
m_{ISO}(M)=\lim_{R\to\infty} m_{ISO}(\Omega(R)).
\ee
Huisken's notion was first introduced with a more complicated definition involving arbitrary exhaustions in  an Oberwolfach report in 2006 \cite{Huisken-Ober-06}. Significantly more details were presented by Huisken in a recorded talk at IAS in 2009.
Huisken's mass was shown to be bounded below by the ADM mass by Miao which
appears in a paper by Fan-Shi-Tam \cite{Fan-Shi-Tam}.  Jauregui-Lee prove the opposite inequality in \cite{Jauregui-Lee-C0-isoper-mass}.   See also Jauregui-Lee-Unger's discussion regarding different versions of the isoperimetric mass in \cite{Jauregui-Lee-Unger}.  They observe that the isoperimetric mass as defined by Huisken is always nonnegative but that proofs that this isoperimetric mass agrees with ADM mass apply the nonnegative scalar curvature.
See also the survey of Benatti-Fogagnolo about isoperimetric notions of mass in lower regularity  \cite{BF-isoper-survey}.   

In addition to the Hawking mass and isoperimetric mass, there are many other quasi-local masses defined for regions which converge to the ADM mass for well selected exhaustions of asymptotically flat manifolds with nonnegative Ricci curvature.  See Liu and Yau's criteria for quasi-local masses in \cite{Liu-Yau-2003,Liu-Yau-quasi} and papers that cite this foundational work.  In particular, there is the Bartnik mass \cite{Bartnik-mass-1986} and the
the Brown-York mass of \cite{Brown-York} studied by Shi-Tam \cite{Shi-Tam}.  These quasilocal masses and their corresponding positivity and rigidity theorems are discussed in the IAS survey about scalar curvature and convergence \cite{Sormani-IAS-survey}.

\section{Examples} \label{sect:examples}

In this section we present examples of sequences of asymptotically flat Riemannian manifolds with nonnegative scalar curvature whose ADM mass approaches zero.   We include graphics and vaguely describe the geometric convergence of these various examples in this section.   Later, in Section~\ref{sect:distances}, when we provide each of the definitions of various notions of geometric convergence, we will refer to these examples and describe exactly how they converge or not towards explicit limit spaces depending on the notion of geometric convergence that is applied.   Here we wish to focus more on the geometry and what kinds of examples can be constructed.   We include open questions about possible additional constructions.

Note that many of these examples are asymptotically flat Riemannian manifolds, $M^3_j$, with nonnegative scalar curvature that contain closed interior minimal surfaces.   Recall that our geometric stability conjecture requires that our manifolds lie in the class, 
$\mathcal{M}$ of Definition~\ref{defn:class}, which is the class defined by Bray in \cite{Bray-Penrose} to prove the Penrose inequality.   This class does not include  $M^3_j$ with closed interior minimal surfaces. Our examples demonstrate why we include this restriction in our conjectures.   We also discuss how cutting along the closed minimal surfaces in $M^3_j$ and taking the connected component, $M'_j$, that includes the asymptotically flat end gives a manifold (called the exterior region) which is in this class, $\mathcal{M}$ of Definition~\ref{defn:class}.

This section has eight subsections.  The first subsection covers Riemannian Schwarzschild space.   The next considers 
scalar flat manifolds with many black holes studied by
Misner, Lichnerowicz, Brill and Lindquist that
they called geometrostatic manifolds.  We then have a subsection on Schoen-Yau tunnels.   The fourth subsection has a review of the spherically symmetric manifolds studied by Lee-Sormani and the construction of a spherical zone in such manifolds.   The fifth subsection covers bubbling and the sixth covers wells.
The sixth subsection covers constructions using tunnels including sewing and scrunching.
The final subsection includes open questions about other possible construction and whether one can construct examples with scrunching without using tunnels.

\subsection{Riemannian Schwarzschild Space}

\begin{figure}[h] 
   \centering
   \includegraphics[width=5in]{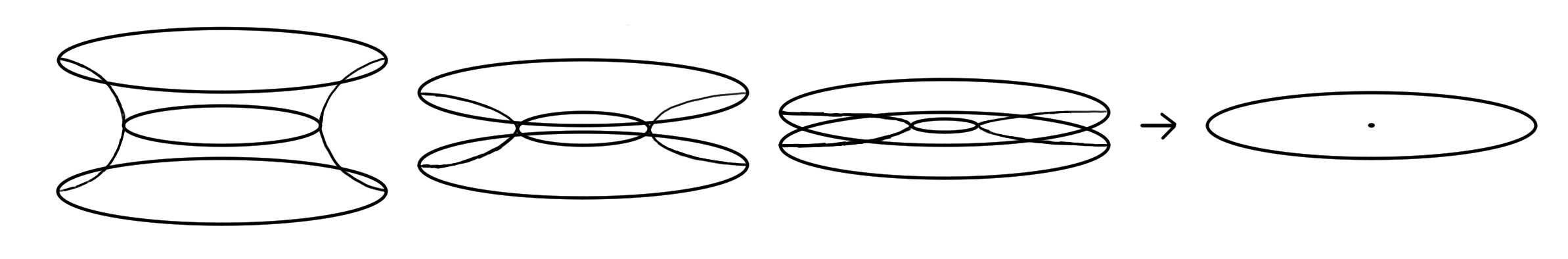} 
   \caption{Riemannian Schwarzschild Spaces with mass converging to zero as in Example~\ref{ex:schwarzschild}.
   }
\label{fig:schwarzschild}
\end{figure}

If one assumes that one has an asymptotically flat Riemannian manifold with zero scalar curvature with a spherically symmetric
metric tensor
then one can solve an ordinary differential equation to conclude that the manifold is either Euclidean space with $m_{ADM}({\mathbb E}^3)=0$ 
or it is Riemannian Schwarzschild space, $M^3_{Sch,m}$,
with $m_{ADM}(M^3_{Sch,m})=m>0$.   In fact
Riemannian Schwarzschild space can be described as a parabola rotated spherically around its indicatrix as depicted in Figure~\ref{fig:schwarzschild}  (c.f. \cite{LeeSormani1} for details). 

\begin{ex} \label{ex:schwarzschild} A sequence of Riemannian Schwarzschild manifolds, $M^3_{Sch,m_j}$, with $m_j \to 0$ as depicted in Figure \ref{fig:schwarzschild} seems to converge geometrically to a pair of Euclidean spaces identified at a point.  Each $M^3_{Sch,m_j}$ has a neck with a spherical minimal surface whose area achieves equality in the Penrose inequality.  One side of this minimal surface defines a manifold $M'_{Sch,m_j}\in \mathcal{M}$ of Definition~\ref{defn:class} with a single end and a boundary that is a minimal surface and no closed interior minimal surfaces.  Examining the sequence as $m_j\to 0$,
it appears that $M'_{Sch,m_j}$ converges geometrically to ${\mathbb E}^3$.
\end{ex}

\subsection{Geometrostatic Manifolds}
Lichnerowicz \cite{Lich-1944},  Misner \cite{Misner-geo} and 
Brill-Lindquist \cite{Brill-Lindquist} studied the following family of complete manifolds with nonnegative scalar curvature with 
 $N+1$ asymptotically flat ends
 that they called geometrostatic manifolds in these papers but should not be confused with more modern definitions of a geometrostatic manifold.

 \begin{ex}\label{ex:geometrostatic}
 A Brill-Lindquist geometrostatic manifold is of the form
 \be
 ({\mathbb R}^3\setminus\{p_1,p_2,...,p_N\},
 g=(\,\chi\,\psi)^2\, g_{{\mathbb E}^3})
 \ee
 where
\be\label{chi-psi}
\chi(x)= 1+\sum_{i=1}^N \frac{\,\alpha_i\,}{\rho_i(x)} \quad\textrm{ and }\quad
\psi(x)=1+\sum_{i=1}^N \frac{\,\beta_i\,}{\rho_i(x)}
\ee
where $\alpha_i$ and $\beta_i$ are positive real numbers and $\rho_i(x)=|x-p_i|$.  Such manifolds have one end as $|x|\to \infty$ whose ADM mass is
\be\label{ADM}
m=m_{N+1}=\sum_{i=1}^N (\alpha_i+\beta_i), 
\ee
and also an end as $x\to p_i$ for
each $i=1,...,N$, whose ADM mass is
\be\label{m_i}
m_i=\alpha_i+\beta_i
+\sum_{k\neq i} \frac{(\beta_i \alpha_k +\beta_k\alpha_i)}{r_{i,k}}
\ee 
where $r_{i,k}=|p_i-p_k|$.   When $N=1$ and
$\alpha_1=\beta_1=m/2$, this is
Riemannian Schwarzschild space of mass $m$ and has a minimal surface at $\rho_1^{-1}(m/2)$.  
\end{ex}

\begin{rmrk}\label{rmrk:Stavrov-Sormani}
Stavrov and Sormani studied the geometric properties of these geometrostatic manifolds
in \cite{Sormani-Stavrov}.  Adding an assumption of a uniform lower bound $\delta>0$
on all the $r_{i,j}$ in a sequence of such manifolds, $(M_j,g_j)$, with ADM mass, $m_{N+1}(M_j)\to 0$, decreasing to $0$, they locate outermost minimal surfaces $\Sigma_i$ in Euclidean annuli about each $p_i$, 
\be
\Sigma_i\subset \rho_i^{-1}(r_i,R_i).
\ee
They thus find the exterior region $M'_j$ that includes the end as $|x|\to \infty$.   In this setting, they easily prove smooth convergence
\be
g_j \to g_{{\mathbb E}^3}
\textrm{ on } {\mathbb R}^3\setminus
\bigcup_{i=1}^N \rho_i^{-1}[0,R_i).
\ee
They control the volumes and areas of regions within the Euclidean annuli sufficiently to prove geometric convergence in the intrinsic flat sense, 
\be
(M'_j,g_j) \Fto ({\mathbb E}^3, g_{{\mathbb E}^3}).
\ee
See Remark~\ref{rmrk:Stavrov-Lakzian}.
\end{rmrk}

See also work by Benko-Stavrov
in \cite{Benko-Stavrov}
and by
Benjamin-McDermott-Stavrov
in \cite{Benjamin-McDermott-Stavrov}
studying sequences of geometrostatic manifolds in a variety of settings.

\subsection{Schoen-Yau Tunnels}

In \cite{Schoen-Yau-tunnels}, Schoen and Yau
proved that one can connect a pair of manifolds, $M_1$ and $M_2$, with positive scalar curvature via a tunnel of positive scalar curvature to create a new manifold of positive scalar curvature, $M_1\# M_2$.  One can imagine this using a pair of three spheres connected by a tunnel similar in shape to the neck in Riemannian Schwarzschild space.   See also Gromov-Lawson \cite{Gromov-Lawson-tunnels}.   Additional important work on gluing and wormholes has been completed by Isenberg-Mazzeo-Pollack \cite{Isenberg-Mazzeo-Pollack}
and Chru\'sciel-Delay
\cite{CD-mapping}.   These tunnels can be applied to create sequences of examples as follows:

\begin{ex}\label{ex:tunnels}
In \cite{Basilio-Dodziuk-Sormani} this tunnel construction is made precise with distance estimates, guaranteeing that arbitrarily small balls can be removed from the standard round, ${\mathbb S}^3$, of constant positive sectional curvature and replaced by arbitrarily short tunnels of positive scalar curvature connecting them.   In further work by Dodziuk appearing in \cite{Dodziuk-tunnels} these tunnels are shown to be able to take prescribed lengths. In work by Sweeney, these tunnels can have arbitrary length and scalar curvature close to that of the spheres if they have sufficiently small width \cite{Sweeney-examples}. 
The volumes of these tunnels can be taken arbitrarily small.
\end{ex}

\subsection{Spherical Zone}

In order to construct examples of sequences of asymptotically flat Riemannian manifolds with nonnegative scalar curvature involving the tunnels of Example~\ref{ex:tunnels},
we will use the following construction.
Details to make it precise can be found in
the Lee-Sormani papers
\cite{LeeSormani1}\cite{LeeSormani2}.

\begin{figure}[h] 
   \centering
   \includegraphics[width=5in]{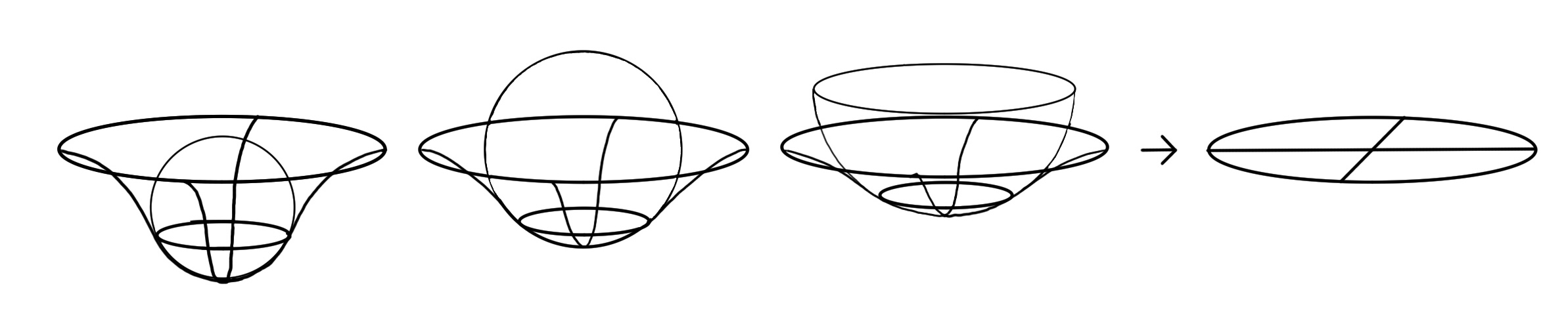} 
   \caption{Here we construct Example~\ref{ex:spherical-zone}.
   }
\label{fig:spherical-zone}
\end{figure}

\begin{ex}\label{ex:spherical-zone}
One can construct a sequence of asymptotically flat spherically symmetric Riemannian manifolds, $M_j^3$, with positive scalar curvature that all contain a spherical zone which is a ball of fixed radius, $R_0>0$, that is isometric to a sphere of constant curvature $K_j \to 0$ as their $m_{ADM}(M_j)\to 0$.   See this construction in Figure~\ref{fig:spherical-zone} which leads to the sequence of manifolds depicted in Figure~\ref{fig:smooth-limit} that converge smoothly to Euclidean space as mass converges to $0$.
\end{ex}

\subsection{Bubbling}

It has been well known that sequences of manifolds with positive scalar curvature can form bubbles.  See for example work of Miao \cite{Miao-PAMS} which builds upon work of Corvino \cite{Corvino-scalar-deformation} and work of
Beig-\'O{} Murchadha
in \cite{BO-trapped}.

\begin{figure}[h] 
   \centering
   \includegraphics[width=5in]{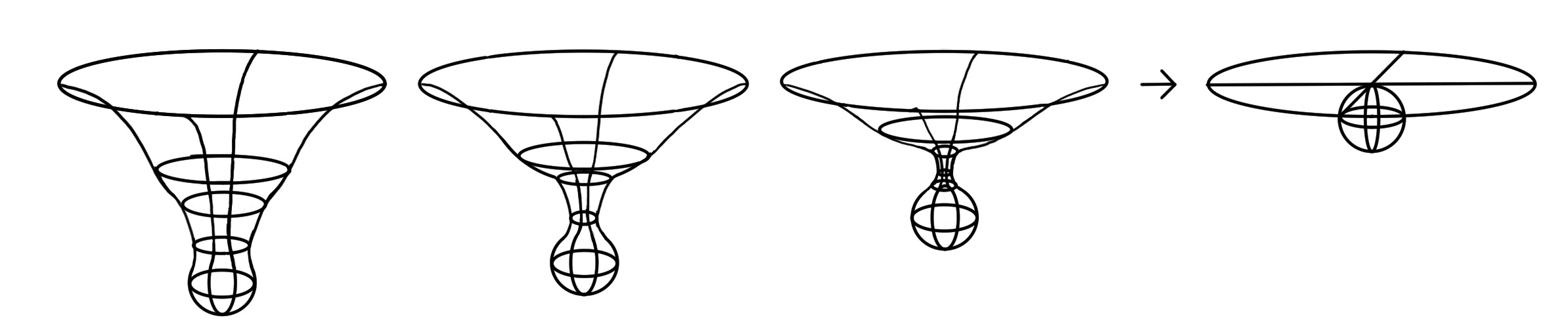}    \caption{Example~\ref{ex:one-bubble}.
}
\label{fig:one-bubble}
\end{figure}

Here is a construction of such a sequence using spherical zones and tunnels:

\begin{ex}\label{ex:one-bubble}
Take $N_j^3$ to be the Riemannian manifolds
with spherical zones constructed in Example~\ref{ex:spherical-zone} and attach each with an increasingly tiny tunnel as described in Example~\ref{ex:tunnels}
to standard spheres, ${\mathbb S}^3$
of fixed size, as depicted in Figure~\ref{fig:one-bubble}.  The resulting sequence of $M_j^3$ with bubbles do not converge smoothly to Euclidean space.   They appear to converge geometrically to Euclidean space with a sphere attached as will be discussed further later.   If we cut along the minimal surfaces in the necks of the tunnels, and study the 
ends' connected components, $M_j'\in \mathcal{M}$ of Definition~\ref{defn:class}, we see that $M_j'$ does appear to converge in some geometric sense to
${\mathbb E}^3$ because the bubbles have been cut off.
\end{ex}

\begin{rmrk}\label{rmrk:Penrose-bubble}
Note that in Example~\ref{ex:one-bubble} we see why the class $\mathcal{M}$ of Definition~\ref{defn:class}  cannot include any interior minimal surfaces
if we hope to prove the Penrose Inequality
as in (\ref{eq:Penrose}).   If one cuts along the equator of the spherical bubble and leaves the sphere in the neck of the bubble, then one has the classical counter example to the Penrose Inequality.  Indeed the area of the equator of the bubble can be of arbitrary size as the mass and area of the minimal sphere in the neck converge to $0$.  To achieve the Penrose inequality, one must restrict to the class, $\mathcal{M}$.
\end{rmrk}

\begin{ex}\label{ex:many-bubbles}
It is also possible to connect increasingly many bubbles to the $N_j^3$ with spherical zones via increasingly many increasingly tiny tunnels to create sequences of $M_j^3$ which have isoperimetric regions, $\Omega_j(R)$,
with volumes diverging to infinity because of all the bubbles.   Luckily one can cut all these bubbles off of the sequence by cutting along minimal surfaces inside the tunnels
and study the 
ends' connected components, $M_j'\in \mathcal{M}$ of Definition~\ref{defn:class}.  Note that by Bray's Penrose Inequality, the total area of $\partial M_j \to 0$.  We see that $M_j'$ does appear to converge in some geometric sense to
${\mathbb E}^3$.
\end{ex}

\begin{rmrk} \label{rmrk:Corvino-bubbles}
In \cite{Corvino-scalar-deformation}, Corvino described a method of locally deforming a manifold with positive scalar curvature.  In particular, he can deform compact regions inside a sequence of Riemannian Schwarzschild spaces with mass decreasing to $0$. 
Miao then applied Corvino's method to
construct sequences with bubbles in \cite{Miao-PAMS}.  See also work of
Beig-\'O{} Murchadha
in \cite{BO-trapped}.   The geometric properties and geometric convergence of these examples have not yet been explored.
\end{rmrk}

\begin{rmrk}\label{rmrk:Anderson-Corvino-Pasqualotto}
Anderson-Corvino-Pasqualotto adapted the
Brill-Lindquist geometrostatic examples
as we described in Example~\ref{ex:geometrostatic},
removing the ends at the $p_i$ and gluing in
smooth regions to create examples that are diffeomorphic to Euclidean space with many bubbles
in \cite{Anderson-Corvino-Pasqualotto}.
\end{rmrk}

\subsection{Wells}

The following construction of wells was well known for some time, popularized in a series of talks by Ilmanen at Columbia.  Some details can be found in Lee-Sormani \cite{LeeSormani1}\cite{LeeSormani2}.
See also a similar construction in Lakzian-Sormani \cite{Lakzian-Sormani}.

\begin{ex}\label{ex:well}
Given a sphere of constant curvature, one can remove a ball of arbitrarily small size and replace it with a spherically symmetric well of positive scalar curvature of arbitrarily small width and arbitrary depth.  If a prescribed bound on depth is provided the volume of this well is less than its boundary's area times this depth.   If no depth bound is provided then the volume can be made arbitrarily large.   These wells contain no minimal surfaces.
\end{ex}

\begin{figure}[h] 
   \centering
   \includegraphics[width=5in]{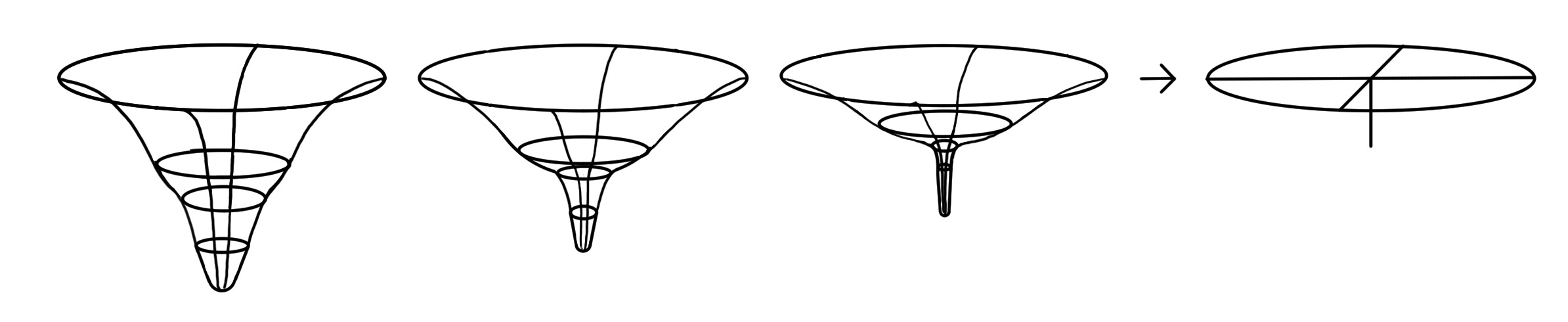} 
   \caption{Example~\ref{ex:one-well}.
   }
\label{fig:one-well}
\end{figure}

\begin{ex}\label{ex:one-well}
Take $N_j^3$ to be the Riemannian manifolds
with spherical zones constructed in Example~\ref{ex:spherical-zone} and attach to each an increasingly thin well of uniform depth $D_0$ as described in Example~\ref{ex:well}.
This creates $M_j^3\in \mathcal{M}$ of Definition~\ref{defn:class} 
 as in Figure~\ref{fig:one-well} which appear to converge geometrically perhaps to Euclidean space attached to a line segment, $[0,D_0]$ at one point.   We will see later that this sequence converges in the intrinsic flat sense to ${\mathbb E}^3$ because the wells volume converges to $0$.
\end{ex}

\begin{figure}[h] 
   \centering
\includegraphics[width=5in]{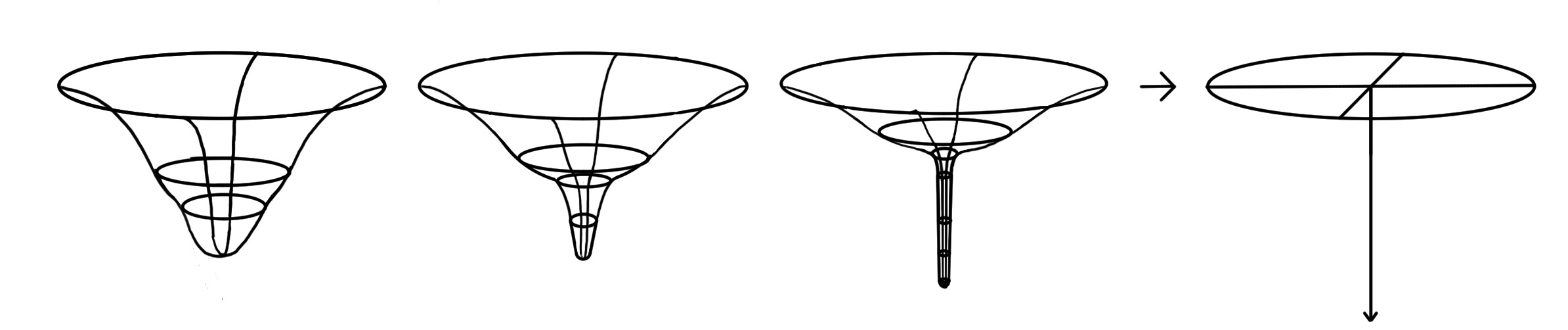} 
\caption{Example~\ref{ex:infinite-well}
}
\label{fig:infinite-well}
\end{figure}

\begin{ex}\label{ex:infinite-well}
Take $N_j^3$ to be the Riemannian manifolds
with spherical zones constructed in Example~\ref{ex:spherical-zone} and attach to each an increasingly thin well of increasing depth $D_j$ as described in Example~\ref{ex:well}.
This creates $M_j^3\in \mathcal{M}$ of Definition~\ref{defn:class} 
as in Figure~\ref{fig:infinite-well} which appear to converge geometrically perhaps to Euclidean space attached to a 
half line at one point.
Depending upon how quickly the depth increases compared to the width this sequence might have $\vol(\Omega_j(R))\to \infty$.   It is one reason we require a diameter bound in Conjecture~\ref{conj:LS}.
\end{ex}

\begin{figure}[h] 
   \centering
\includegraphics[width=5in]{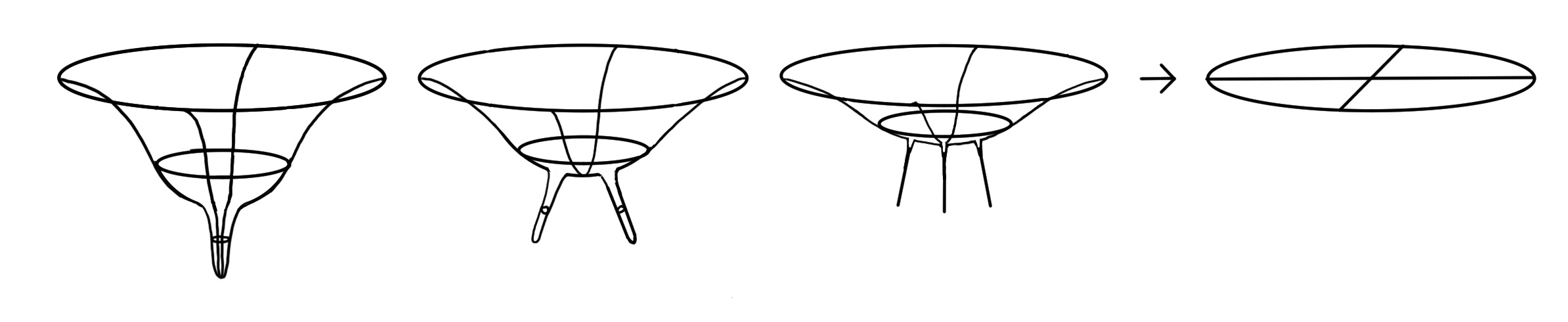} 
\caption{Example~\ref{ex:many-wells}.
   }
\label{fig:many-wells}
\end{figure}

\begin{ex}\label{ex:many-wells}
Take $N_j^3$ to be the Riemannian manifolds
with spherical zones constructed in Example~\ref{ex:spherical-zone} and attach to each increasingly many increasingly thin wells of uniform depth $D_0$ as described in Example~\ref{ex:well}.
This creates $M_j^3\in \mathcal{M}$ of Definition~\ref{defn:class} 
as in Figure~\ref{fig:many-wells} which we will later explain have no GH limit but converges in the $\mathcal{F}$ sense
to Euclidean space because the total filling volume of these wells converges to $0$.
\end{ex}

\begin{rmrk} \label{rmrk:Corvino-wells}
In \cite{Corvino-scalar-deformation}, Corvino describes a method of locally deforming a manifold with positive scalar curvature.   This method can also be applied to create wells.  The geometric properties of Corvino's deformations are not yet well explored.  They provide a wealth of additional examples of sequences of manifolds which are isometric to Schwarzschild spaces outside of compact regions which may or may not have bubblings, wells, or other geometric properties.
\end{rmrk}

\subsection{Sewing}

\begin{figure}[h] 
   \centering
\includegraphics[width=5in]{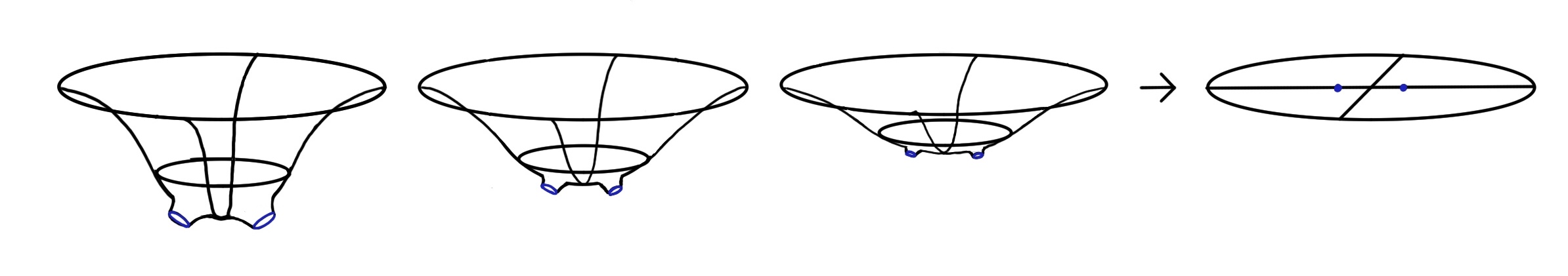} 
\caption{Example~\ref{ex:one-tunnel}.
   }
\label{fig:one-tunnel}
\end{figure}

Here we will see a collection of examples
constructed using the tunnels described in Example~\ref{ex:tunnels}.  This time rather than running the tunnels between two different manifolds, we will run the tunnel between two balls in the same manifold.   We begin with a simple example:

\begin{ex}\label{ex:one-tunnel}
Take $N_j^3$ to be the Riemannian manifolds
with spherical zones constructed in Example~\ref{ex:spherical-zone} and remove pairs of increasingly small balls from the spherical zones that are a fixed distance $D_0>0$ apart.
Replace each pair of balls with an increasing tiny tunnel as in Example~\ref{ex:tunnels}.
This creates $M_j^3$
as in Figure~\ref{fig:one-tunnel} which appears to converge geometrically perhaps to Euclidean space with two points identified as a single point in the limit space.
If we cut along the minimal surfaces inside each tunnel, we create
$M_j'\in \mathcal{M}$ of Definition~\ref{defn:class} which converge to ${\mathbb E}^3$ because there is no longer a short cut running through the tunnels.
\end{ex}

In \cite{Basilio-Dodziuk-Sormani}, Basilio-Dodziuk-Sormani described a method called "sewing along a curve" in which they edit in increasingly small increasingly short tunnels along a curve in a three-sphere to obtain a sequence of manifolds with positive scalar curvature that converge to a sphere with a pulled string.
Such pulled string spaces were first described by Burago while working with Ivanov and Sormani on early versions of \cite{Burago-Ivanov-area}.  A sphere with a pulled string is a sphere with an entire curve identified to a single point.   See the precise description of the placement of the tunnels to achieve this limit in \cite{Basilio-Dodziuk-Sormani}.  It is not difficult to adapt their sewing along a curve to construct the following example.

\begin{figure}[h] 
   \centering
  \includegraphics[width=5in]{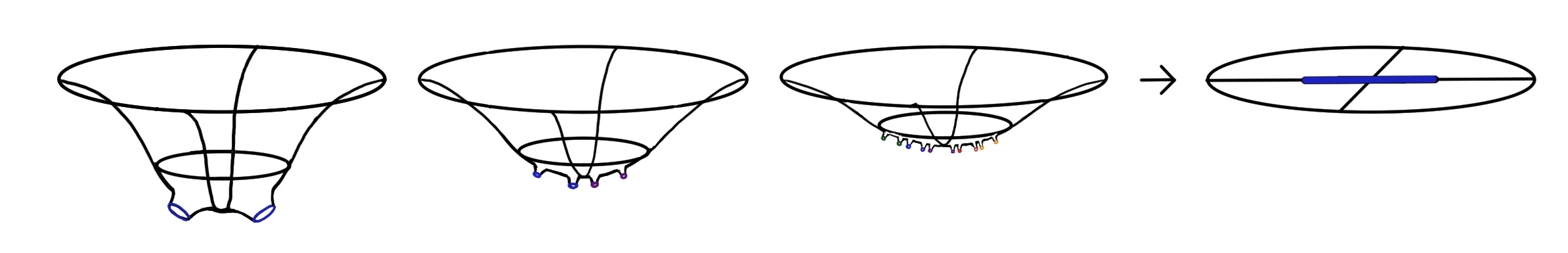}    \caption{Example~\ref{ex:sewing-curve}
   constructed by sewing along a curve as in work of Basilio-Dodziuk-Sormani.
}
\label{fig:sewing-curve}
\end{figure}

\begin{ex}\label{ex:sewing-curve}
Starting with the sequence of manifolds, $N^3_j$ with spherical zones as in Example~\ref{ex:spherical-zone}, one can take a sequence of curves, $C_j$, of length $L_0$, in the spherical zones converging to a curve, $C_\infty$, of length $L_0$ in Euclidean space.   One may then construct a sequence of  increasingly many increasingly thin and short tunnels sewn along the curves $C_j$ as in Basilio-Dodziuk-Sormani \cite{Basilio-Dodziuk-Sormani} to define manifolds $M_j^3$ as depicted in Figure~\ref{fig:sewing-curve}.  The sequence of $M^3_j$ converge geometrically to  Euclidean space with the curve, $C_\infty$, identified to a point. 
If we cut along the minimal surfaces in the tunnels, we obtain $M_j'\in \mathcal{M}$ which converge to Euclidean space because the short cuts through the tunnels have been removed.
\end{ex}

In \cite{Basilio-Sormani}, Basilio-Sormani also described a method of "sewing" a region to a point in which they edit in an increasingly small increasingly dense network of tunnels across a region in a three-sphere to obtain a sequence of manifolds with positive scalar curvature that converge to a three sphere with the entire region identified to a point.   It is not difficult to adapt their sewing along a curve to construct the following example.

\begin{ex}\label{ex:sewing-region}
Starting with the sequence of manifolds, $N_j^3$, with spherical zones as in Example~\ref{ex:spherical-zone}, one can take a sequence of regions in the spherical zones converging to a region in Euclidean space.   One may then construct a network of increasingly many increasingly thin and short tunnels sewn into the regions as in Basilio-Sormani \cite{Basilio-Sormani}.   The limit space would then be Euclidean space with a region identified to a point.
\end{ex}

\begin{rmrk}
Other examples of manifolds with positive scalar curvature involving tunnels that are increasingly thin but long can be found in the work of Basilio-Kazaras-Sormani, Sweeney, and Krandel-Sweeney \cite{BKS-no-geod}
\cite{Sweeney-examples}
\cite{Krandel-Sweeney-long}.
It is possible that with carefully selected tunnels
similar to those in the work of Krandel-Sweeney applied to the spherical zone of Example~\ref{ex:spherical-zone}, one could obtain a limit space which is a smooth Riemannian manifold that is flat Euclidean space outside of a ball yet has shorter distances between points than
flat Euclidean space.
\end{rmrk}

\begin{rmrk}\label{rmrk:no-tunnels}
It should be noted that sequences with tunnels like
Example~\ref{ex:one-tunnel},
Example~\ref{ex:sewing-curve}
and Example~\ref{ex:sewing-region}
do not contradict 
Conjecture~\ref{conj:LS}
despite the fact that they approach limits 
which are not Euclidean space.  This is because they have closed minimal surfaces inside the tunnels, so they are not in the class $\mathcal{M}$ of Definition~\ref{defn:class}.  
\end{rmrk}

\begin{rmrk}\label{rmrk:cut-tunnels}
When one cuts the minimal surfaces in Example~\ref{ex:one-tunnel},
Example~\ref{ex:sewing-curve}
and Example~\ref{ex:sewing-region}
to create exterior manifolds that lie in the class $\mathcal{M}$ of Definition~\ref{defn:class}, those exteriors do converge to Euclidean space because there are no longer short cuts running through the tunnels that lead to shortened distances in
the limit spaces.
\end{rmrk}

\begin{rmrk}
Additional tunnel constructions appear in the work of
Isenberg-Mazzeo-Pollack \cite{Isenberg-Mazzeo-Pollack}
and Chru\'sciel-Delay
\cite{CD-mapping}. These papers have an analytic approach and their geometric properties have not yet been well studied.  It would be intriguing to explore the geometry of these important classes of manifolds.
\end{rmrk}

\subsection{Scrunching Possibilities}
\label{sect:scrunching}

In \cite{Basilio-Sormani}, Basilio-Sormani define a concept called {\em scrunching}.   In their paper they take a fixed manifold, $M^3$, and a subset $A_0\subset M^3$. 
They say that a sequence of manifolds,
\be
N_j^3= (M^3 \setminus A_{\delta_j})\disjointunion A'_{\delta_j}, 
\ee
is said to {\em ``scrunch $A_0$ down to a point''}
if $A_{\delta_j}$ is the tubular neighborhood, $T_{\delta_j}(A_0)$ and
$A'_{\delta_j}$ is a smooth scrunched replacement that
satisfies:
\be\label{sewn-curve-tubular'}
\vol(A_{\delta_j}')\le \vol(A_{\delta_j})(1+\epsilon_j), 
\ee
\be\label{sewn-curve-vol'}
\vol(N_j^3) \le \vol(M^3) (1+\epsilon_j)
\ee
and
\be \label{sewn-curve-diam'}
\diam(A_{\delta_j}')\le H_j
\ee
where 
\be\label{eq:BS-to-0}
\epsilon_j \to 0
\qquad H_j \to 0, \qquad 2\delta_j<H_j.
\ee
Note that this definition of scrunching does not mention scalar curvature, nor does it describe how the convergence of diameter to zero is achieved.

Basilio-Sormani prove in \cite{Basilio-Sormani}
that under the conditions above, including (\ref{sewn-curve-tubular'})-(\ref{sewn-curve-diam'}), the scrunched sequence
$N_j^3$ converges to $M^3$ with the scrunched region, $A_0$, identified to a point in the
$GH$, $mm$ and $\mathcal{F}$ sense. 
Note that their theorem about scrunching does not mention scalar curvature, and is purely a geometric convergence theorem that can be applied in many settings.   
In that same paper, Basilio-Sormani apply this scrunching theorem to prove the convergence of their examples with positive scalar curvature 
similar to Example~\ref{ex:sewing-region} which are constructed using the sewing method with  increasingly many increasingly small tunnels.   In the remarks at the end of this subsection we discuss possibilities of constructing scrunching regions without using the sewing method.

\begin{figure}[h] 
   \centering  \includegraphics[width=5in]{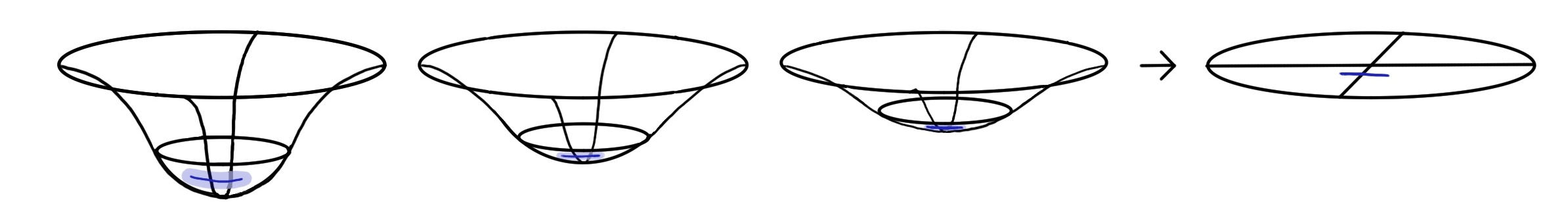} 
   \caption{A sequence of manifolds with $ADM$ mass converging to $0$ that have sets (depicted as blue curves) which ``scrunch'' to a point so that the limit space is Euclidean space with a set identified to a point.   
}
\label{fig:scrunching}
\end{figure}

Here we will slightly adapt the Basilio-Sormani definition of scrunching, allowing ourselves to consider smoothly converging, $N_j \to {\mathbb E}^3$, as in Example~\ref{ex:spherical-zone} with spherical zones
depicted in Figure~\ref{fig:spherical-zone}.  We take
$A_0\in {\mathbb E}^3$ and $A_j\subset M_j$ smoothly converging to
$A_0$ as in Figure~\ref{fig:scrunching}.   Then we
say that a sequence of manifolds,
\be
M_j^3= (N_j^3 \setminus A_{j,\delta_j})\disjointunion A'_{j,\delta_j}, 
\ee
is said to {\em ``scrunch $A_0$ down to a point''}
if $A_{j,\delta_j}$ is the tubular neighborhood, $T_{\delta_j}(A_j)\subset N^3_j$ and
$A'_{j,\delta_j}$ is a smooth scrunched replacement that
satisfies (\ref{sewn-curve-tubular'})-(\ref{sewn-curve-diam'}).   It should be easy to prove (imitating the work by Basilio-Sormani in \cite{Basilio-Sormani}) that 
the new scrunched sequence
$N_j^3$ converges $GH$, $mm$ and $\mathcal{F}$ sense
to $M^3$ with the scrunched region, $A_0$, identified to a point.  This new theorem can then be applied using various techniques, either with sewing of tunnels or without tunnels, to prove new examples related to Geometric Stability of the Zero Mass Theorem.   See the remarks below.

\begin{rmrk}\label{rmrk:scrunching}
If one can construct a sequence of $M_j^3$ as in Figure~\ref{fig:scrunching} with nonnegative scalar curvature and ADM mass converging to $0$ and no closed interior minimal surfaces 
which scrunch a region to a point so that their limit is Euclidean space with a region identified to a point, then
one would have a counter example to Conjecture~\ref{conj:LS}. Note that the Basilio-Sormani method of sewing a region to obtain scrunching has tunnels and many closed interior minimal surfaces.
Nobody has yet been able to produce scrunching without tunnels and minimal surfaces in three dimensions with nonnegative scalar curvature. 
\end{rmrk}

\begin{rmrk}\label{rmrk:Corvino-scrunching}
Intuitively, one might imagine that a region with scrunching would have to contain significant amounts of positive scalar curvature perhaps concentrating along a curve that one hopes to pull to a point in the limit.  Perhaps one might try to complete such a construction using the methods of Corvino in \cite{Corvino-scalar-deformation} by prescribing scalar curvature inside the spherical zones in the sequence in Example~\ref{ex:spherical-zone}.   It would also be of interest to prove that Corvino's method cannot be applied to shrink distances between points in this way.
\end{rmrk}

\begin{rmrk}\label{rmrk:Mantoulidis-Schoen}
Another intriguing construction is that of Mantoulidis-Schoen in \cite{Mantoulidis-Schoen-Bartnik}.  They deform a compact region near the spherical horizon in the exterior Riemannian Schwarzshild space to create necks that end at new horizons which are not just standard spheres.   This work has been extended by Cabrera Pacheco, Cederbaum, McCormick and Miao
in \cite{CCMM-Bartnik} to create regions with horizons of a variety of shapes.
It appears that the necks in these deformations can possibly have large diameters, however, if one can control the depths of these perturbations then perhaps one can use these methods to create scrunched regions without an horizon.
Cederbaum has a team exploring depth bounds in these constructions.
\end{rmrk}

\begin{rmrk}\label{rmrk:new-notions}
If a scrunching counter example 
as described in Remark~\ref{rmrk:scrunching} were constructed, one would need to explore alternative notions of Geometric Convergence like those described in the work of Lee-Naber-Neumayer \cite{LNN-dp} and Dong-Song \cite{Dong-Song-Stability} as we will describe below.
One might also consider a geometric notion of convergence
based on areas alone as first explored by Sormani with Burago and Ivanov with some initial progress published in \cite{Burago-Ivanov-area}.
\end{rmrk}

\section{Geometric Notions of Convergence}
\label{sect:distances}

In this section we review a variety of geometric notions of convergence, including Gromov-Hausdorff, metric measure, and intrinsic flat convergence. We immediately apply each notion to the various examples of sequences of Riemannian manifolds satisfying the hypotheses of Conjecture~\ref{conj:LS} on the Geometric Stability of the Schoen-Yau Zero Mass Rigidity Theorem.   We encourage the reader to consider developing other possible geometric notions of convergence that might fit this setting as well.

Most geometric notions of convergence are defined by using definite distances between pairs of compact Riemannian manifolds, $d_{dist}(M_1,M_2)$ such that
\be
d_{dist}(M_1,M_2)=0 \quad \iff
\quad M_1 \textrm{ is isometric to } M_2.
\ee
Then one defines convergence, 
\be
M_j \distto M_\infty \iff d_{dist}(M_j,M_\infty)\to 0.
\ee

To study convergence to noncompact limit spaces like Euclidean space,
$M_j \distto M_\infty$ one studies convergence of 
precompact domains, $\Omega_j(R)\subset M_j$,
from a
canonical exhaustions of $M_j$ to 
precompact domains, $\Omega_{\infty}(R)\subset M_\infty$,
from a
canonical exhaustions of $M_\infty$.  When the domains are balls around fixed points $p_j\in M_j$, this is called pointed convergence and the limit space may depend upon the choice of points.
However, one may also consider convergence of isoperimetric regions or some other natural exhaustion of the manifolds.   The various papers mentioned within all consider different canonical exhaustions so be sure to read the papers to find the precise statements.

To define a distance between a pair of compact Riemannian manifolds,  $(\Omega_j,g_j)$, one first converts them into metric spaces,
$(\Omega_j, d_{g_j})$, using the Riemannian distance as reviewed in (\ref{eq:Riem-dist-1})-(\ref{eq:Riem-dist-2}).   Note that metric spaces do not have coordinate charts, however one can define metric differentials as in
work of Korevaar-Schoen \cite{Korevaar-Schoen}.
See also work of 
Cheeger, Heinonen, Kirchheim, Shanmugalingam, and
others (cf. \cite{Cheeger-diff}\cite{Heinonen-nonsmooth} 
\cite{Kirchheim}\cite{Shanmugalingam-Newtonian}).

After converting the Riemannian manifolds into metric spaces, one considers a pair of maps into an arbitrary complete metric space, $(Z,d_Z)$,
\be \label{eq:dist-pres-pair}
\varphi_j: (\Omega_j, d_j)\to (Z, d_Z),
\ee
that are distance preserving as in (\ref{eq:dist-pres}).   One can then measure the distance between the images of these distance preserving maps,
\be 
d_{ext-dist}^Z(\varphi_1(\Omega_1),\varphi_2(\Omega_2)),
\ee
using various extrinsic notions of distances that depend upon how these images lie within $Z$.   This approach was first introduced by 
Edwards \cite{Edwards} and then Gromov \cite{Gromov-metric}.

To achieve an intrinsic notion of distance between the original pair of manifolds, one takes an infimum over all possible $(Z,d_Z)$ and all possible distance preserving maps into $Z$:
\be \label{eq:inf-ext-dist}
d_{dist}(\Omega_1,\Omega_2)=\inf d_{ext-dist}^Z(\varphi_1(\Omega_1), \varphi_2(\Omega_2)).
\ee
For the Gromov-Hausdorff distance, $d_{GH}$, one uses the Hausdorff distance, $d_H^Z$ in $Z$ (see Section~\ref{sect:GH}),
for the metric measure distance, $d_{mm}$, one may use the Wasserstein distance in $Z$, $d^Z_W$, (see Section~\ref{sect:mm}), and for the intrinsic flat distance, $d_{\mathcal F}$, one uses the Federer-Fleming flat distance, $d^Z_F$, in $Z$ (see Section~\ref{sect:F}).   The later two notions need additional structure on the metric spaces which will be described within these sections.   

\begin{rmrk}\label{rmrk:definite}
We encourage the development of additional notions of distances between Riemannian manifolds using other extrinsic distances within complete metric spaces, $Z$, to control other aspects of the geometry of the space.
To study rigidity and geometric stability, one needs the notion of distance between Riemannian manifolds to be definite in the sense that
\be
d_{dist}(\Omega_1,\Omega_2)=0 \quad \implies \quad
\Omega_1 \textrm{ is isometric to } \Omega_2.
\ee
This can sometimes be achieved 
when the extrinsic distance in
$Z$ satisfies
\be
d_{ext-dist}^Z(K_1,K_2)=0 \quad
\implies \quad K_1=K_2
\ee
however one must be careful to show the infimum
in the definition of $d_{dist}$ is achieved.
If additional structure is desired beyond just an
isometry, $\Psi: \Omega_1\to \Omega_2$, then that additional
structure must be controlled by  $d_{ext-dist}^Z$
and must interact well with the infimum.
\end{rmrk}

\begin{rmrk}\label{rmrk:another}
Another way to develop a stronger notion of convergence is to combine multiple notions.   For example one may take the sum of the GH and mm distances to ensure measured Gromov-Hausdorff convergence.  The volume preserving intrinsic flat convergence introduced by Portegies in \cite{Portegies-F-evalues} is defined by
\be
d_{\mathcal{VF}}(\Omega_1,\Omega_2)=
d_{\mathcal F}(\Omega_1,\Omega_2) + |\vol(\Omega_1)-\vol(\Omega_2)|
\ee
See also Sormani's older survey articles \cite{Sormani-survey-scalar} 
and \cite{Sormani-IAS-survey}. 
\end{rmrk}

\subsection{Gromov-Hausdorff Convergence}
\label{sect:GH}

\begin{figure}[h] 
   \centering
\includegraphics[width=5in]{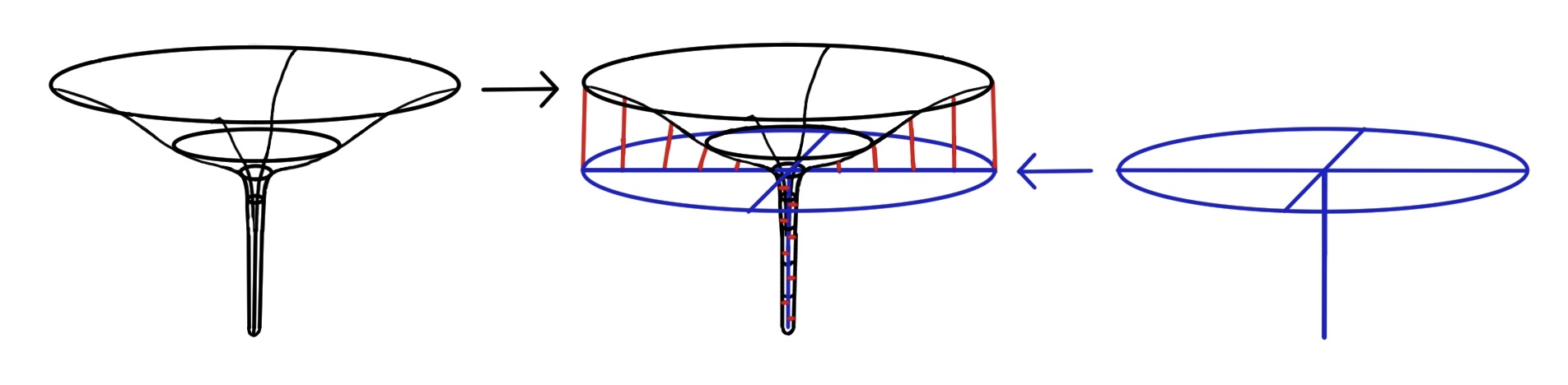} 
   \caption{Here we estimate the GH distance between $\Omega_j(R)\subset M_j$ on the left with a well and  the GH limit, $\Omega_\infty(R)\subset M_\infty$, 
   of Example~\ref{ex:one-well} which is a Euclidean ball with a line segment attached on the right.
They are embedded together in $Z$ in the center where the red vertical and horizontal segments are estimating the Hausdorff distance between them in $Z$.  Notice how the segment in $M_\infty$
must be close in length to the depth of the well in $M_j$ to achieve a small GH distance between them.
   }
\label{fig:GH-well}
\end{figure}

Gromov-Hausdorff convergence was first introduced by 
Edwards \cite{Edwards} and then reintroduced by Gromov \cite{Gromov-metric} who proved a corresponding compactness theorem.
The  Gromov-Hausdorff distance between a pair of compact Riemannian manifolds, $(\Omega_j,g_j)$, is defined by converting them into metric spaces, $(\Omega_j, d_j)$,
as in (\ref{eq:Riem-dist-1})-(\ref{eq:Riem-dist-2})
and
taking an infimum over all complete metric spaces, $(Z,d_Z)$, and over all distance preserving maps, 
\be
\varphi_j: (\Omega_j, d_j)\to (Z, d_Z),
\ee
of the Hausdorff distance between the images,
\be
d_{GH}(\Omega_1,\Omega_2)=\inf\, d_H^Z(\varphi_1(\Omega_1),\varphi_2(\Omega_2))
\ee
where the Hausdorff distance between compact subsets, $K_1,K_2\subset Z$, is
defined by
\be\label{eq:Hausdorff}
d^Z_H(K_1,K_2)=\inf 
\left\{ R\ \middle\vert \begin{array}{l}
 \forall \, x\in K_1 \, \exists \, y_{x,R}\in K_2
\textrm{ s.t. } d_Z(x,y_{x,R})<R \\
    \forall \, y\in K_2 \,\exists \, x_{y,R}\in K_1
\textrm{ s.t. } d_Z(y,x_{y,R})<R
  \end{array}\right\}.
\ee
See Figure~\ref{fig:GH-well} where $K_1$ is the flat disk with a line segment attached in blue, $K_2$ is the surface with the well in black, and the red segments are estimating the Hausdorff distance between them, thus providing an estimate on the Gromov-Hausdorff distance between the manifolds on the left and on the right.

Gromov proved this is a definite notion:
\be
d_{GH}(\Omega_1,\Omega_2)=0 \quad \iff \quad \Omega_1
\textrm{ is isometric to } \Omega_2.
\ee
He also proved a compactness theorem which can be applied to sequences of Riemannian manifolds with a uniform lower bound on their Ricci curvature \cite{Gromov-metric}.

Example~\ref{ex:schwarzschild} depicted in Figure~\ref{fig:schwarzschild} of Riemannian Schwarzschild manifolds, $M_j$, with $\mass_{ADM}(M_j)\to 0$ converges in the pointed GH sense to a pair of Euclidean spaces jointed at a point.   If one considers only the exterior Schwarzschild (outside the horizon), they converge in the pointed GH sense to Euclidean space.

Example~\ref{ex:spherical-zone} depicted in Figure~\ref{fig:spherical-zone} of asymptotically flat Riemannian manifolds, $M_j$, with spherical zones of constant curvature and $m_{ADM}(M_j)\to 0$ converges in the pointed GH sense to Euclidean space.   In fact, this example converges smoothly and in every sense discussed in this article.

Example~\ref{ex:one-well}  of 
isoperimetric regions, $\Omega_j(R)$, in asymptotically flat Riemannian manifolds, $M_j$, with single wells of depth $D$ and $m_{ADM}(M_j)\to 0$ converges in the pointed GH sense to a ball $B_0(R)$ in Euclidean space with an interval $[0,D]$ attached to it as depicted in Figure~\ref{fig:one-well}.   In Figure~\ref{fig:GH-well} we see how one term in the sequence is close to the GH limit with the segment.  

Example~\ref{ex:infinite-well}  of asymptotically flat Riemannian manifolds, $M_j$, with single wells of depth $D_j \to \infty$ and $m_{ADM}(M_j)\to 0$ converges in the pointed GH sense to Euclidean space with a half line $[0,\infty)$ attached to it as depicted in Figure~\ref{fig:infinite-well}.   

Example~\ref{ex:many-wells} depicted in Figure~\ref{fig:many-wells} of asymptotically flat Riemannian manifolds, $M_j$, with increasingly many wells of depth $D$ and $m_{ADM}(M_j)\to 0$ has no limit in the pointed GH sense since the sequence is not equicompact.

Example~\ref{ex:one-bubble} of asymptotically flat Riemannian manifolds, $M_j$, with single bubbles of radius $R_j\to R$ and $m_{ADM}(M_j)\to 0$ converges in the pointed GH sense to Euclidean space with a sphere of radius $R$ attached  as depicted in Figure~\ref{fig:one-bubble}.        

Example~\ref{ex:one-tunnel} depicted in Figure~\ref{fig:one-tunnel} of asymptotically flat Riemannian manifolds, $M_j$, with single short tunnels between points a distance $R$ apart and $m_{ADM}(M_j)\to 0$ converges in the pointed GH sense to Euclidean space with two points identified.    


Example~\ref{ex:sewing-curve} depicted in Figure~\ref{fig:sewing-curve} of asymptotically flat Riemannian manifolds, $M_j$, with many short tunnels sewn along a curve and $m_{ADM}(M_j)\to 0$ converges in the pointed GH sense to Euclidean space with a curve identified to a point.  The tunnels have created short cuts along the entire curve so that distances between points along the curve converge to $0$.   The detailed proof is completed by Basilio-Dodziuk-Sormani \cite{Basilio-Dodziuk-Sormani}.

Example~\ref{ex:sewing-region} of asymptotically flat Riemannian manifolds, $M_j$, with many short tunnels in a network covering a region and $m_{ADM}(M_j)\to 0$ converges in the pointed GH sense to Euclidean space with the region identified to a point.   This is proven in detail by Basilio-Sormani in \cite{Basilio-Sormani} along with a detailed proof of the GH convergence of scrunched regions like those described in Remark~\ref{rmrk:scrunching}.

It is important to note that the Gromov-Hausdorff distance only controls distances, not volumes or areas or boundaries.   While there are many beautiful results about Gromov-Hausdorff limits which do control these additional geometric properties, including work of Cheeger-Colding, Petersen-Wei and others CITES
this is done with additional curvature bounds which imply stronger convergence properties like, for example, metric measure convergence.

\begin{rmrk}\label{rmrk:harmonic-maps}
In \cite{KKL}, Kazaras-Khuri-Lee study the Gromov-Hausdorff stability of zero mass rigidity theorem by assuming a uniform lower bound on the Ricci curvature of the $M_j$.
They apply the ADM mass bound of
Bray-Kazaras-Khuri-Stern in \cite{BKKS} which is proven building on work of Stern in \cite{Stern-scalar} that controls harmonic maps using scalar curvature. Note that all the examples of sequences with increasingly thin tunnels and wells fail to satisfy the uniform lower bound on Ricci curvature.  Nevertheless the estimates in these papers have lead to more general advances by Dong and
Dong-Song in \cite{Dong-some-stability}\cite{Dong-Song-Stability}
which will be discussed in Remark~\ref{rmrk:Dong-Song}.   See also
work of Hirsch-Miao-Tsang
in \cite{HMT-CAG}.
\end{rmrk}

\subsection{Metric Measure Convergence}
\label{sect:mm}

The notion of metric measure convergence of Riemannian manifolds was first introduced by Fukaya in \cite{Fukaya-collapsing} and studied further in work of Cheeger-Colding including \cite{ChCo-PartI}.  Sturm defined the metric measure or D-distance between metric measure spaces in \cite{Sturm-2006-I}.   Sturm and Lott-Villani developed the notion of $CD$ metric measure spaces with generalized nonnegative Ricci curvature in \cite{Sturm-2006-I} and \cite{Lott-Villani-09}.   Ambrosio-Gigli-Savar\'e developed the notion of $RCD$ metric measure spaces with a stronger notion of generalized nonnegative Ricci curvature
in \cite{AGS} building upon work by Gigli that was posted online years ago but only recently published in \cite{Gigli-splitting}.

In \cite{Fukaya-collapsing}, Fukaya proposed converting compact Riemannian manifolds into metric measure spaces,
$(\Omega_j, d_{g_j}, \mu_{g_j})$, with a Borel probability measure
\be \label{eq:prob}
\mu_j(A)=\vol_{g_j}(A)/\vol_{g_j}(\Omega_j),
\ee
which is called the renormalized volume measure.
More generally one can define metric measure spaces,
$(M,d,\mu)$, where $(M,d)$ is a complete separable metric space and $\mu$ is a Borel probability measures.   

There are a few ways to define metric measure convergence which do not exactly agree with one another.  We will use the simplest notion that uses distance preserving maps and pushforward measures. Given a Lipschitz map, $\phi:X\to Z$, into a 
metric spaces, $Z$, one can define
the push forward of the measure,
\be 
\phi_*\mu(A)=\mu(\phi^{-1}(A)).
\ee
Sturm defined the D-distance or metric measure distance between metric measure spaces in  \cite{Sturm-2006-I} by 
taking an infimum over all complete metric spaces, $(Z,d_Z)$, and over all distance preserving maps, 
\be
\varphi_j: (X_j, d_j)\to (Z, d_Z),
\ee
of the $L^2$-Wasserstein distance between the pushforwards of their measures,
\be\label{eq:mm}
d_{{mm}}(\Omega_1,\Omega_2)=\inf\, d_W^Z(\varphi_{1*}\mu_1,\varphi_{2*}\mu_2).
\ee
The Wasserstein distance between the
measures is 
\be
d_W (\nu_1,\nu_2) := \inf_{q\in P(M\times M)}
\left(\int_{M\times M} d^2(x,y) dq(x,y)\right)^
{1/2}
\ee
where the infimum is taken over all probability
measures $q$ on the product space $M\times M$ whose image measures under the projections
$\pi_1(x,y)=x$ and $\pi_2(x,y)=y$
satisfy
\be
\pi_{1*}q=\nu_1 \textrm{ and } \pi_{2*}q=\nu_2.
\ee
The Wasserstein distance essentially measures the
$L^2$ cost of transporting the measure $\nu_1$ to $\nu_2$. It is achieved by an optimal transport plan which describes optimal paths between the points where the measures have positive density.  See Figure~\ref{fig:GH-well} where 
we estimate the metric measure distance between
the flat disk on the right and the surface with the well on the left.  In the center we see the sets of positive density of the measures that have been pushed forward into a metric space, $Z$.
The red segments are describing the optimal paths
in $Z$,  needed to estimate on the Wasserstein distance between the push-forward measures.

\begin{figure}[h] 
   \centering
\includegraphics[width=5in]{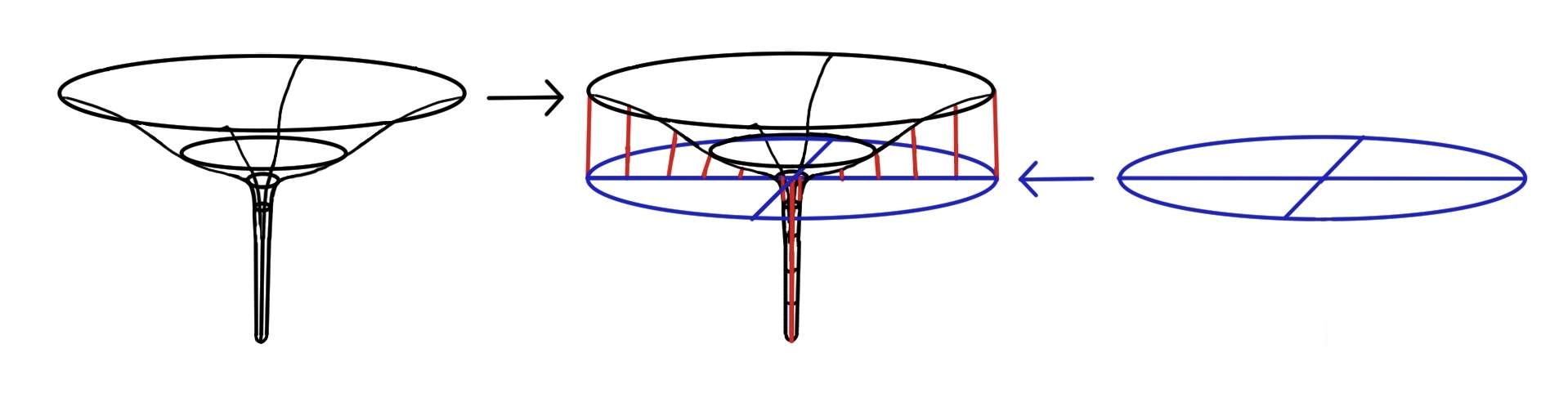} 
   \caption{Here we estimate the Wasserstein distance between $\Omega_j(R)\subset M_j$ on the left with a well and the metric measure limit, $\Omega_\infty(R)\subset M_\infty$, 
   of Example~\ref{ex:one-well} which is a Euclidean ball, $B_0(R)\subset {\mathbb E}^3$. 
   Their renormalized volume measures
   as in (\ref{eq:prob}) are pushed forward into a common metric space, $Z$, at the center.
   Notice how the image of the well in $\Omega_j(R)$
has very small volume and so it costs
very little to transport it upward along the red
optimal path up to the image of the disk.   So we achieve a small mm distance between them.}
\label{fig:mm-well}
\end{figure}

Sturm proved this is a definite notion:
$$
d_{mm}(\Omega_1,\Omega_2)=0 \iff
\exists \textrm{ an isometry }F:\Omega_1\to \Omega_2
\textrm{ s.t. } F_*\mu_1=\mu_2.
$$
He also proved a compactness theorem which can be applied to sequences of Riemannian manifolds with uniform lower bounds on their Ricci curvature \cite{Sturm-2006-I}.

Example~\ref{ex:schwarzschild} depicted in Figure~\ref{fig:schwarzschild} of Riemannian Schwarzschild manifolds, $M_j$, with $m_{ADM}(M_j)\to 0$ converges in the mm sense to a pair of Euclidean spaces jointed at a point.   If one considers only the exterior Schwarzschild (outside the horizon), they converge in the mm sense to Euclidean space.   Example~\ref{ex:spherical-zone} depicted in Figure~\ref{fig:spherical-zone} of asymptotically flat Riemannian manifolds, $M_j$, with spherical zones of constant curvature and $m_{ADM}(M_j)\to 0$ converges in the mm sense to Euclidean space. 

Example~\ref{ex:one-well}  of 
isoperimetric regions, $\Omega_j(R)$, in asymptotically flat Riemannian manifolds, $M_j$, with single wells of depth $D$ and $m_{ADM}(M_j)\to 0$ converges in the mm sense to a ball $B_0(R)$ in Euclidean space.  This is not the GH limit.  The wells disappear in the $mm$ limit because their volumes converge to $0$.  See Figure~\ref{fig:mm-well}   

Example~\ref{ex:one-bubble} of asymptotically flat Riemannian manifolds, $M_j$, with single bubbles of radius $R_j\to R$ and $m_{ADM}(M_j)\to 0$ converges in the pointed mm sense to Euclidean space with a sphere of radius $R$ attached  as depicted in Figure~\ref{fig:one-bubble}.        

Example~\ref{ex:one-tunnel} depicted in Figure~\ref{fig:one-tunnel} of asymptotically flat Riemannian manifolds, $M_j$, with single short tunnels between points a distance $R$ apart and $m_{ADM}(M_j)\to 0$ converges in the mm sense to Euclidean space with two points identified.    

Example~\ref{ex:sewing-curve} depicted in Figure~\ref{fig:sewing-curve} of asymptotically flat Riemannian manifolds, $M_j$, with many short tunnels sewn along a curve and $m_{ADM}(M_j)\to 0$ converges in the mm sense to Euclidean space with a curve identified to a point.  The tunnels have created short cuts along the entire curve so that distances between points along the curve converge to $0$.   The detailed proof is completed by Basilio-Dodziuk-Sormani \cite{Basilio-Dodziuk-Sormani}.

Example~\ref{ex:sewing-region} of asymptotically flat Riemannian manifolds, $M_j$, with many short tunnels in a network covering a region and $m_{ADM}(M_j)\to 0$ converges in the mm sense to Euclidean space with the region identified to a point only if the volume of the region is $0$.   This is proven in detail by Basilio-Sormani in \cite{Basilio-Sormani} along with a detailed proof of the mm convergence of scrunched regions pf zero volume like those described in Remark~\ref{rmrk:scrunching}.

It is important to note that the metric measure distance only controls distances and volumes but not areas.     In addition,
metric measure convergence does not naturally control boundaries of sets and does not interact particularly well with isoperimetric regions.  It is perhaps best to use metric convergence combined with another notion of convergence to achieve geometric stability of the zero mass rigidity theorem. 

\subsection{Sormani-Wenger Intrinsic Flat Convergence}
\label{sect:F}

Sormani-Wenger introduced intrinsic flat convergence in \cite{SW-JDG} building upon work of Ambrosio-Kirchheim in \cite{AK}.  
To find the intrinsic flat distance between
a given pair of compact oriented $m$-dimensional Riemannian manifolds with a boundary of finite area, one first converts them into integral current spaces,
$(\Omega_j, d_{g_j}, [[\Omega_j]])$, with an integral current structure, $[[\Omega_j]]$, defined by its action on $m-forms$,
\be\label{eq:current-structure}
[[\Omega_j]](\omega)=\int_{\Omega_j} \omega
\ee
and a boundary, $\partial [[\Omega_j]]$, defined 
by its action on $(m-1)-forms$,
\be
\partial[[\Omega_j]](\sigma)=
[[\Omega_j]](d\sigma)=\int_{\Omega_j} d\sigma
=
\int_{\partial \Omega_j}\sigma.
\ee
More general integral current spaces, $(X,d,T)$, are defined by Sormani-Wenger in \cite{SW-JDG} using Ambrosio-Kirchheim's currents on complete metric spaces \cite{AK}.
A current, $T$, in a complete metric space, $Z$, acts on tuples of Lipschitz functions, $\pi_i:Z\to {\mathbb R}$, 
with finite Ambrosio-Kircheim mass, $\mass_{AK}(T)$.  When
$T=[[M]]$, we have
\be
T(\pi_0, \pi_1,...,\pi_m)=\int_M \pi_0\, d\pi_1\wedge\cdots\wedge d\pi_m
\ee
and
\be
\mass_{AK}([[M]])=\vol_{g}(M)
\textrm{ and } \mass_{AK}(\partial[[M]])=\vol_{g}(\partial M).
\ee
In general, an integral current has boundary
defined by,
\be
\partial T(\sigma_1,...,\sigma_{m-1})=
T(1,\sigma_1,...,\sigma_{m-1}),
\ee
and pushforwards by Lipschitz maps, $\phi:Z\to W$, are defined by
\be
\phi_\#T(\pi_0,\pi_1,...,\pi_m)=T(\pi_0\circ \phi,
\pi_1\circ \phi,...,\pi_m\circ \phi).
\ee
This notion is essential if one is interested in keeping track of boundaries, slicing by Lipschitz functions, volumes, and areas, and rectifiable structures \cite{SW-JDG}.

\begin{figure}[h] 
   \centering
\includegraphics[width=5in]{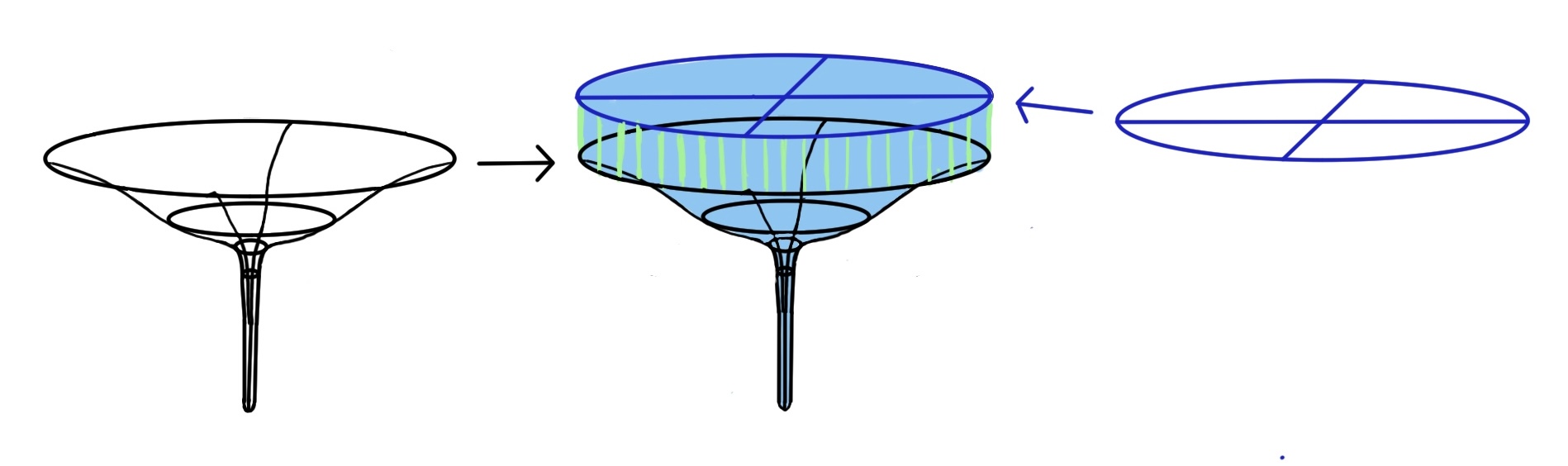} 
   \caption{Here we see $\Omega_j(R)\subset M_j$ on the left with a well and on the right the $\mathcal{F}$ limit, $B_0(R)\subset {\mathbb E}^3$, 
   of Example~\ref{ex:one-well}.
They are embedded together in $Z$ in the center where the blue region is the filling volume between them
with volume, $\mass_{AK}(B)$, and
the green cylinder has volume, $\mass_{AK}(A)$.  Notice how the 
depth of the well does not contribute significantly to this $\mathcal{F}$ distance between the spaces as long as it is very thin and has small filling volume.
   }
\label{fig:F-well}
\end{figure}

The intrinsic flat distance between a pair of compact oriented Riemannian manifolds, $(\Omega_j,g_j)$, is defined by converting them into integral current spaces, $(\Omega_j, d_j, [[\Omega_j]])$,
as in (\ref{eq:Riem-dist-1})-(\ref{eq:Riem-dist-2}) and
(\ref{eq:current-structure})
and
taking an infimum over all complete metric spaces, $(Z,d_Z)$, and over all distance preserving maps, 
\be
\varphi_j: (\Omega_j, d_j)\to (Z, d_Z),
\ee
of the Flat distance between the pushforwards of their current structures,
\be\label{eq:intrinsic-flat}
d_{{\mathcal F}}(\Omega_1,\Omega_2)=\inf\, d_F^Z(\varphi_{1\#}[[\Omega_1]],\varphi_{2\#}[[\Omega_2]])
\ee
where the flat distance between integral currents, $T_1,T_2$,
in $Z$ is
defined by
\be\label{eq:flat}
d^Z_F(T_1,T_2)=\inf 
\left\{ \mass(A)+\mass(B)\,:\,
A+\partial B=T_1-T_2\right\}.
\ee
This flat distance is essentially
a weighted filling volume, $\mass(B)$,
between the $T_j=\varphi_{j\#}[[\Omega_j]]$
 with an extra current, $A$, 
 added in so that we have Stokes Theorem, 
 \be
 \partial B=T_1-T_2-A.
\ee
See the center of Figure~\ref{fig:F-well} where $T_1$ on top is the flat disk in blue, $T_2$ on the bottom is the black surface with a well, $A$ is the cylinder in green, and $B$ is the blue interior, so that
$\mass_{AK}(A)+\mass_{AK}(B)$
provides an estimate on the $\mathcal{F}$ distance between the manifolds on the left and on the right.

In \cite{SW-JDG}, Sormani and Wenger proved the intrinsic flat distance is a definite notion:
\be
d_{\mathcal{F}}(\Omega_1,\Omega_2)=0 \quad \iff \quad \exists  
\Psi: \Omega_1\to \Omega_2
\ee
where $\Psi$ is an orientation preserving isometry.
A slightly revised definition of the intrinsic flat distance which rescales well under scalings of distances in the pair of spaces  introduced by LeFloch-Sormani in \cite{LeFloch-Sormani} is also definite.

Wenger proved a compactness theorem
which states that
any sequence of
compact oriented $\Omega_jj$ such that
\be \label{eq:Wenger}
\vol(\Omega_j)\le V
\textrm{ and }
\area(\partial \Omega_j) \le A
\textrm{ and }
\diam(\Omega_j)\le D
\ee
has a subsequence which converges in the $\mathcal{F}$ sense to an integral current space \cite{Wenger-compactness}.  Note that this
compactness theorem extends the compactness theorems of Federer-Fleming and of Ambrosio-Kirchheim \cite{FF}\cite{AK} but only requires a 
diameter bound rather than a uniform compact set.
For complete asymptotically flat manifolds, $M_j$, we will consider the intrinsic flat convergence of the uniform isoperimetric regions, $\Omega_j(R)$. 

Example~\ref{ex:one-well} depicted in Figure~\ref{fig:one-well} of 
isoperimetric regions, $\Omega_j(R)$, in asymptotically flat Riemannian manifolds, $M_j$, with single wells of depth $D$ and $m_{ADM}(M_j)\to 0$ converges in the $\mathcal{F}$ sense to a ball $B_0(R)$ in Euclidean space .   In Figure~\ref{fig:F-well} we see how one term in the sequence is close to the ball in Euclidean space.  A detailed proof including a construction of $Z$ and an explicit estimate on the filling volume is provided in work of Lee-Sormani \cite{LeeSormani1}.   In fact, they proved intrinsic flat stability of the positive mass theorem in the spherically symmetric setting using such explicit constructions.   

Example~\ref{ex:many-wells} depicted in Figure~\ref{fig:many-wells} of asymptotically flat Riemannian manifolds, $M_j$, with increasingly many wells of depth $D$ and $m_{ADM}(M_j)\to 0$ has no limit in the pointed GH sense but does converge in the intrinsic flat sense to Euclidean space.  This can be seen by filling in all the wells and showing the volume of the total filling converges to $0$.  

Example~\ref{ex:one-tunnel}, Example~\ref{ex:sewing-curve}, and
Example~\ref{ex:sewing-region} created with increasingly short tunnels, as well as Example~\ref{ex:one-bubble} with 
a bubble all converge in the ${\mathcal F}$ and GH sense to the same limit spaces, as in their figures.   This is because the increasingly small tunnels have increasingly small filling volumes.
The details can be found in Dodziuk-Basilio-Sormani and Basilio-Sormani 
\cite{Basilio-Dodziuk-Sormani}
\cite{Basilio-Sormani}.  Note that Remark~\ref{rmrk:no-tunnels}, Remark~\ref{rmrk:cut-tunnels}, and Remark~\ref{rmrk:scrunching} all relate to both $\mathcal{F}$ and $GH$ convergence.

Sormani-Wenger proved in \cite{SW-JDG} that
\be
\Omega_j \Fto \Omega_\infty
\ee
implies
\be
\liminf_{j\to\infty}
\vol(\Omega_j) \ge 
\vol(\Omega_\infty)
\ee
when $\Omega_\infty$ is smooth, that
\be
\partial \Omega_j \to \partial \Omega_\infty
\ee
and that
\be
\liminf_{j\to\infty}
\area(\partial \Omega_j) \ge 
\area(\partial \Omega_\infty).
\ee
In \cite{Sormani-AA},
Sormani proved Arzela-Ascoli type theorems about uniformly Lipschitz functions,
\be
F_j: \Omega_j \to N
\ee
having subsequences
which converge to a
Lipschitz function
\be
F_\infty:\Omega_\infty\to N.
\ee
Portegies-Sormani have proven many additional properties of intrinsic flat convergence in \cite{Sormani-Portegies-prop} using
the methods of Ambrosio-Kirchheim \cite{AK}.

\begin{rmrk}\label{rmrk:GMT-approach}
Huang-Lee-Sormani studied geometric stability of the zero mass rigidity theorem for manifolds which are graphs in Euclidean space satisfying various technical properties in \cite{HLSa}.  They considered regions $\Omega_j\subset M_j'$ where
$\mass(M_j)\to 0$ and applied Wenger's Compactness theorem to obtain a subsequential $\mathcal{F}$ limit space $\Omega_\infty$.   Applying earlier work 
of Huang-Lee
\cite{HL} and Lam \cite{Lam-graph}, they could show
$\vol(\Omega_j)$ converges to the Euclidean volume.  They then applied semicontinuity of volumes, the Arzela-Ascoli Theorem, and control on the boundaries claiming to demonstrate that
$\Omega_\infty$ must be isometric to Euclidean Space.   This final claim was made rigorous by Del-Nin and Perales in \cite{DNP}
using advanced techniques from the Geometric Measure Theory on metric spaces developed by Ambrosio and Kirchheim in \cite{AK,AK2,Kirchheim}.   See also work of Basso, Creutz and Soultanis \cite{BCS}.
 Further work applying this approach in more general settings is in progress by Perales and her collaborators.   
\end{rmrk}

\subsection{Volume Preserving Intrinsic Flat Convergence}
\label{sect:VF}

Portegies introduced the notion of volume preserving intrinsic flat 
$\mathcal{VF}$
convergence in \cite{Portegies-F-evalues}.   We can say
$\Omega_j\VFto \Omega_\infty$ iff
\be
d_{\mathcal VF}(\Omega_j,\Omega_\infty)
=d_{\mathcal F}(\Omega_j,\Omega_\infty)+
|\vol(\Omega_j)-\vol(\Omega_\infty)|\to 0.
\ee

Every example presented in this paper has the same $\mathcal{VF}$ limit as its $\mathcal F$ limit because the volumes are all converging.  This is one reason why the zero mass stability conjecture originally stated by Lee-Sormani in \cite{LeeSormani1} for $\mathcal F$ convergence is now stated for $\mathcal VF$ convergence in Conjecture
\ref{conj:LS} within
this paper and also in the surveys
\cite{Sormani-survey-scalar}
\cite{Sormani-IAS-survey}.

There are strong consequences resulting when a sequence of manifolds converges in the volume-preserving intrinsic flat sense.
Portegies proved that eigenvalues of the Laplacian with Neumann boundary conditions semi-converge in \cite{Portegies-F-evalues}.   Jauregui-Lee proved that ADM mass semi-converges in \cite{Jauregui-Lee-VF-isoper-mass} by studying good competitors for isoperimetric regions.
Jauregui-Perales-Portegies proved semi-continuity of capacity in \cite{Jauregui-Perales-Portegies-capacity}.
These powerful consequences of $\mathcal{VF}$ convergence justify why we strive to achieve volume-preserving intrinsic flat convergence in our conjectures.   

For those wishing to test the geometric stability conjecture as stated in 
Conjecture~\ref{conj:LS} with the class of manifolds, $\mathcal{M}$ of Definition~\ref{defn:class}, it is worth checking if one can perhaps directly prove the following consequences of this conjecture: 

\begin{conj}\label{conj:vol}
Under the hypotheses of Conjecture~\ref{conj:LS},
\be\vol(\Omega_j(R))\to \vol(B_0(R))=4\pi R^3/3.
\ee
See Shi's volume inequality for one direction
\cite{Shi-isoper}.
\end{conj}

\begin{conj}\label{conj:evalue}
Under the hypotheses of Conjecture~\ref{conj:LS},
the eigenvalues of the Laplacian of  
$\Omega_j(R)$ semi-converge to the 
eigenvalues of the Laplacian of
the Euclidean ball, $B_0(R)$ as described by Portegies in \cite{Portegies-F-evalues}.
\end{conj}

\begin{conj}\label{conj:cap}
Under the hypotheses of Conjecture~\ref{conj:LS},
the capacity of $\Omega_j$ 
as in work of Jauregui-Perales-Portegies in \cite{Jauregui-Perales-Portegies-capacity}
semi-converges 
to the capacity of
the Euclidean ball, $B_0(R)$.
\end{conj}

A counter example to any of these three conjectures may lead to a counter example to Conjecture~\ref{conj:LS}.  Such a
counter example might lead one to restrict the class of manifolds, $\mathcal{M}$, of Definition~\ref{defn:class} in this conjecture.  We are deliberately vague in our conjectures and strongly encourage the reader to consult the work of Portegies, Jauregui, Lee, and Perales cited above before working on these conjectures.

\subsection{Smooth, Lipschitz, and VADB convergence}

All the Geometric Notions of Convergence we've described above allow one to consider limits of sequences of manifolds which are not diffeomorphic to one another.  However there are special settings where one considers only sequences of manifolds which are already diffeomorphic to the proposed limit space.   In these settings one will often study a fixed manifold, perhaps simply a ball in ${\mathbb R}^3$ with a sequence of metric tensors, $g_j$, on that ball which have nonnegative scalar curvature.   

If $g_j\to g_\infty$ in the $C^0$ or smoother sense then it is known that if the limit $g_\infty$ is smooth it must have nonnegative scalar curvature as well.  This is proven by Gromov in \cite{Gromov-Dirac} and by Bamler in \cite{Bamler}.  

If one only knows that $g_j\to g_0$ uniformly,
one can conclude that the resulting metric spaces converge in the Lipschitz sense,
\be\label{eq:Lip-conv}
d_{Lip}((\Omega_j,d_j),(\Omega_\infty,d_\infty))
=\inf \log\left(\Lip(\Psi_j)+
\Lip(\Psi_j^{-1})\right)\to 0
\ee
where the infimum is taken over all bi-Lipschitz maps 
$\Psi_j:\Omega_j\to \Omega_\infty$.

When $g_\infty$ is not $C^2$ smooth, one cannot define nonnegative scalar curvature in the usual way. However, 
Gromov in \cite{Gromov-Dirac} and Burkhardt-Guim in \cite{Burkhardt-Guim-pointwise} have proposed lower regularity versions of scalar curvature that can be applied in such a setting.   See \cite{Sormani-IAS-survey} and \cite{Gromov-four} for more details.

Gromov proved that Lipschitz convergence as in (\ref{eq:Lip-conv}) implies GH convergence in \cite{Gromov-metric} and
Sormani-Wenger proved Lipschitz convergence implies $\mathcal{F}$ convergence in \cite{SW-JDG}.   Since Lipschitz convergence implies volume convergence (with the appropriate notion of volume on the limit space) one has $mm$ and ${\mathcal VF}$ convergence as well.

It should be noted that sequences of manifolds with increasingly thin wells as in Examples~\ref{ex:one-well} and ~\ref{ex:many-wells} do not converge in this Lipschitz sense.  An alternative notion of convergence has been introduced by Allen-Perales-Sormani in \cite{APS-VADB} called Volume Above Distance Below (VADB) convergence which works well for  examples with wells. 
The original definition of VADB convergence applies to a manifold without boundary, $(M,g_j)$,
and we say 
\be
g_j\VADBto g_\infty
\ee
iff
\be \label{eq:VADB}
\vol_{g_j}(M)\to
\vol_{g_\infty}(M)
\qquad
g_j \ge g_\infty
\qquad
\diam_{g_j}(M)\le D.
\ee
Allen-Perales-Sormani
proved VADB convergence implies $\mathcal{VF}$ convergence in \cite{APS-VADB}.
Allen-Perales and Allen-Bryden extended this notion and this theorem to convergence of manifolds with boundary in \cite{Allen-Perales-bndry-VADB} and \cite{AB-VADBV2}.   

\begin{rmrk}\label{rmrk:VADB}
VADB convergence methods have been applied by Allen-Perales and by Huang-Lee-Perales to prove special cases of $\mathcal{VF}$-geometric stability of the zero mass rigidity theorem 
in the graph setting \cite{Allen-Perales-bndry-VADB}, \cite{HLP}. 
It has the potential to be applied to a much larger class of spaces.   The real challenge is to prove that one can bound the distances from below by Euclidean distances.  
\end{rmrk}

\subsection{Lebesgue and Sobolev vs Geometric Convergence}
\label{sect:Sobolev}
\label{sect:LNN}

In many settings assumptions on scalar curvature enable one to prove the metric tensors converge in a Sobolev sense to a limit metric tensor $g_j$.  
For example, Allen has proven 
Lebesgue and Sobolev stability of the zero mass rigidity theorem for regions covered by smooth inverse mean curvature flow in $M_j'\subset \mathcal{M}$ 
using work of Huisken-Ilmanen \cite{Huisken-Ilmanen} in 
\cite{Allen-IMCF-23}.
Bryden has proven Sobolev stability of the zero mass rigidity theorem for axisymmetric manifolds away from the axis
in \cite{Bryden-stability} applying work of Brill \cite{Brill-axi} and Chrusciel \cite{Chrusciel-axi}.  In these papers the metric tensors are shown to converge in the Sobolev sense to a Euclidean metric.

However Sobolov convergence of metric tensors does not imply geometric convergence without imposing additional controls.
Examples and Theorems by Allen, Bryden, and Sormani in
\cite{AS-contrasting, AS-relating}, \cite{Allen-Lp-GH}
\cite{Allen-Bryden-Sobolev-Ineq}
explore sequences of metric tensors $g_j$ which converge in Lebesgue and Sobolev senses to limit spaces.  These papers do not assume nonnegative scalar curvature.   It should be noted there are examples, particularly in \cite{AS-contrasting} which have scrunching and pulled strings in the limit spaces even without tunnels.   So one would need to apply nonnegative scalar curvature and perhaps also asymptotic flatness and mass to zero to avoid the construction of such counter examples to the geometric stability of the zero mass rigidity theorem.

In \cite{LNN-dp}\cite{LNN-survey}, Lee-Naber-Neumayer
convert Riemannian manifolds $(M,g)$ into metric spaces using 
\be
d_{g,p}(x,y)=\sup\left\{|f(x)-f(y)|\,:\,\int_M|\nabla f|^p\, d\vol_g\,\le \, 1\,\right\}
\ee
rather than the Riemannian distance, $d_g$.
This $d_{g,p}$ distance controls $W^{1,p}$ Sobolev spaces well and has the property that
\be
\lim_{p\to \infty} d_{g,p}(x,y)=d_g(x,y).
\ee
They then define $d_p$ convergence of Riemannian manifolds using the measured GH convergence of the metric spaces and apply this to study manifolds with nonnegative scalar curvature.   

\begin{rmrk}
\label{rmrk:dpF}
One might also consider the intrinsic flat distances between Riemannian manifolds that have been converted into metric spaces using the Lee-Naber-Neumayer $d_{p,g}$ distance between points in the manifolds.   
\end{rmrk}

\begin{rmrk}
\label{rmrk:dplim}
It is not yet known if one can prove a Lee-Naber-Neumayer style $d_p$ stability for the zero mass rigidity theorem.  In fact the $d_p$ limits of the various examples presented in this paper are unknown.
\end{rmrk}

\subsection{Converging Away from Bad Sets}
\label{sect:smooth-away}

The easiest way to prove $GH$ and $\mathcal{F}$ convergence using methods from Geometric Analysis is to prove that regions within the manifolds have smooth or at least $C^0$ convergence and then control distances, or volumes and areas of the sets where there is no smooth convergence.

\begin{rmrk}\label{rmrk:Lakzian}
In \cite{Lakzian-Sormani}, Lakzian-Sormani prove a theorem in which they take a pair of Riemannian manifolds $\Omega_j$ with regions, $U_j\subset \Omega_j$, where the metric tensors on the $U_j$ are close within $\epsilon_j$
and where
$\lambda_j$ estimating the difference between intrinsic and extrinsic distances,
\be\label{eq:Lakzian-lambda}
\lambda_j=\max_{p,q\in U_j}|d_{\Omega_j}(p,q)-d_{U_j}(p,q)|
\ee
is small.   They conclude that
$
d_{GH}(\Omega_1,\Omega_2)
$
is small depending on
$\epsilon_j$,
$\lambda_j$,
and
\be
h_j=d_{GH}(\Omega_j,U_j)
\ee
which is essentially the
radius $R_j$ needed for any point in $\Omega_j$ to reach a point in $U_j$.
Applying this estimate to Example~\ref{ex:one-well} taking $U_j\subset \Omega_j$ to be a region smoothly close to a Euclidean annulus, then $h_j$ is the depth of the well.
Lakzian-Sormani also conclude that
$
d_{\mathcal{F}}(\Omega_1,\Omega_2)
$
is small depending on
$\epsilon_j$,
$\lambda_j$,
and
\be
V_j=\vol_j(\Omega_j\setminus U_j)
\ee
and
\be
A_j=\area_j(\partial U_j).
\ee
Applying this estimate to Example~\ref{ex:one-well} taking $U_j\subset \Omega_j$ to be a region smoothly close to a Euclidean annulus, then $A_j$ is the area of the opening of the well and $V_j$
is the volume of the well.
\end{rmrk}

\begin{rmrk}
\label{rmrk:Stavrov-Lakzian}
Sormani-Stavrov applied the Lakzian-Sormani theorem to prove $\mathcal{VF}$ stability of the Zero Mass Rigidity Theorem in the setting of Brill-Lindquist Geometrostatic Manifolds in \cite{Sormani-Stavrov}.  Recall that a Brill-Lindquist Geometrostatic manifold, $M$, is a conformally flat manifold with nonnegative scalar curvature on ${\mathbb R}^3\setminus \{p_1,\ldots,p_N\}$ with $N+1$ asymptotically flat ends as described in Example~\ref{ex:geometrostatic}. To complete their proof, Sormani-Stavrov first had to locate the minimal surfaces to find exterior regions, $M'\in \mathcal{M}$, satisfying the hypothesis of the Penrose Inequality and the hypothesis of the Geometric Stability of the Zero Mass Theorem.  They cut out regions with well-controlled volumes and boundary areas around these minimal surfaces as described in Remark~\ref{rmrk:Stavrov-Sormani}.   This allowed them to control the $\epsilon_j$ and $\delta_j$
 needed to apply Lakzian-Sormani's theorem described in Remark~\ref{rmrk:Lakzian}.
They carefully compare lengths of curves in
the spaces to control $\lambda_j$ of (\ref{eq:Lakzian-lambda}) using the fact that $g_j \ge g_{{\mathbb E}^3}$ and the removed balls are easy to traverse around.
They thus prove the sequence converges in the $\mathcal{F}$ sense to a ball in Euclidean space.   
\end{rmrk}

\begin{rmrk}
\label{rmrk:Dong-Song}
Dong-Song have also studied the Geometric Stability of the Zero Mass Rigidity Theorem by removing regions
in \cite{Dong-Song-Stability}.  They apply the ADM mass bound of
Bray-Kazaras-Khuri-Stern in \cite{BKKS} which is proven building on work of Stern in \cite{Stern-scalar} that controls harmonic maps using scalar curvature.
They do not cut along minimal surfaces, so the manifolds they consider can have bubbles like those in Examples~\ref{ex:one-bubble} and~\ref{ex:many-bubbles} and can also have tunnels like those in 
Examples~\ref{ex:one-tunnel} and~\ref{ex:sewing-curve}.   Dong-Song prove ptGH convergence of asym flat regions $U_j$ so they control $\epsilon_j$, 
and they prove $\area(\partial U_j)\to 0$, but cannot prove $\vol(M_j\setminus U_j)\to 0$ because they do not exclude the examples with bubbling.  They cannot control $\lambda_j$ of (\ref{eq:Lakzian-lambda})
because they do not exclude the sewing examples with shortcuts in distances.
\end{rmrk}

\begin{rmrk}
\label{rmrk:Dong-Song-convergence}
It is natural to ask whether Dong-Song have defined a new geometric notion of convergence.  Is there a distance between Riemannian manifolds, $d_{DS}(M_1,M_2)$, that is definite and converges to $0$ for sequences that converge in the sense described by Dong-Song in \cite{Dong-Song-Stability}.   Does this notion of convergence have a compactness theorem?   What are the properties of this notion of convergence?
\end{rmrk}

\begin{rmrk}
\label{rmrk:Dong-Song-F}
It is also natural to ask whether the results in Dong-Song can be made strong enough to achieve $\mathcal{F}$ convergence by requiring that the manifolds lie in $\mathcal{M}$ of Definition~\ref{defn:class} as in the hypotheses of the Penrose Inequality and Conjecture~\ref{conj:LS}.   Song, who is an expert in $\mathcal{F}$ convergence, has stated that such an achievement would require a significant advance beyond the work they have done.   It is unclear how to use outermost minimal surfaces to avoid the possibility of scrunching as described in
Remark~\ref{rmrk:scrunching} where $\lambda_j$ of (\ref{eq:Lakzian-lambda}) is not converging to $0$. 
\end{rmrk}

We would like to close with an apology to everyone whose work has not been mentioned in this paper.   We have not attempted to present all the alternative proofs of the Mass Comparison Theorem nor the lower regularity results nor the analytical results even in dimension three.   We have not touched upon this subject in higher dimensions nor included any discussion of asymptotically hyperbolic manifolds.  We also apologize that we have not provided detailed statements of the results that we have presented.  We hope that we have succeeded in providing a strong intuitive understanding of the challenges one faces trying to prove geometry stability of Zero Mass Rigidity Theorem even in dimension three and that this will provide some insight into related questions as well.

\bibliographystyle{plain}
\bibliography{2025}

\end{document}